\pdfoutput=1

\documentclass[11pt,a4paper, twoside]{amsart}
\usepackage{geometry}
\usepackage[utf8]{inputenc}
\usepackage[english]{babel}
\usepackage[T1]{fontenc}
\usepackage{lmodern}

\usepackage{amsmath}
\usepackage{amssymb}
\usepackage{amsthm}

\usepackage{mathtools}
\usepackage{stmaryrd}

\usepackage{wrapfig}
\usepackage{subcaption}

\usepackage{enumitem}
\usepackage{verbatim}
\usepackage{wrapfig}
\usepackage{makecell}
\usepackage[dvipsnames]{xcolor}

\usepackage{hyperref}
\hypersetup{colorlinks,linkcolor={blue!75!black},citecolor={blue!75!black},urlcolor={blue!75!black}}

\usepackage[backend=biber,maxnames=50]{biblatex}
\AtEveryBibitem{%
	\clearfield{number}}

\usepackage{tikz}

\DeclareMathOperator*{\id}{id}
\DeclareMathOperator*{\im}{im}

\DeclareMathOperator*{\supp}{supp}

\DeclareMathOperator*{\linh}{span}

\DeclareMathOperator*{\conv}{conv}

\DeclareMathOperator*{\aff}{aff}

\DeclareMathOperator*{\bd}{bd}

\newcommand{\MID}{\,\middle|\,}
\DeclarePairedDelimiter\scpr{\langle}{\rangle}

\newcommand{\norm}[1]{\|#1\|}

\newcommand{\R}{\mathbb{R}}

\newcommand{\N}{\mathbb{N}}

\mathtoolsset{showonlyrefs}

\theoremstyle{theorem}

\newtheorem{theorem}{Theorem}[section] %
\newtheorem{corollary}[theorem]{Corollary}
\newtheorem{proposition}[theorem]{Proposition}
\newtheorem{lemma}[theorem]{Lemma}

\theoremstyle{definition}

\newtheorem{definition}[theorem]{Definition}

\theoremstyle{remark}
\newtheorem{remark}[theorem]{Remark}

\numberwithin{equation}{section}

\DeclareMathOperator*{\vertices}{vert}

\DeclareMathOperator*{\relint}{relint}
\DeclareMathOperator*{\vol}{vol}

\newcommand{\volt}[1]{{\textstyle\vol_{#1}}}

\newcommand{\transp}{{\mkern-1.5mu\mathsf{T}}}

\newcommand{\E}{\mathbb{E}}

\newcommand{\sphere}{\mathbb{S}^{n-1}}
\newcommand{\Prob}{\mathbb{P}}
\newcommand{\dint}{\,\mathrm{d}}

\let\emptyset\varnothing

\title{Isotropic constants and regular polytopes}
\author{Christian Kipp}

\address{Technische Universität Berlin, Institut für Mathematik, Sekr.~MA4-1, Straße des 17.~Juni 136, 10623 Berlin, Germany}
\email{kipp@math.tu-berlin.de}
\begin{document}
	
	\begin{abstract}
		We discuss first-order optimality conditions for the isotropic constant and combine them with RS-movements to obtain structural information about polytopal maximizers. Strengthening a result by Rademacher, it is shown that a polytopal local maximizer with a simplicial vertex must be a simplex. A similar statement is shown for a centrally symmetric local maximizer with a simplicial vertex: it has to be a cross-polytope. Moreover, we show that a zonotope that maximizes the isotropic constant and that has a cubical zone must be a cube. Finally, we consider the class of zonotopes with at most $n+1$ generators and determine the extremals in this class.
	\end{abstract}
	
	\maketitle
	
	\setcounter{tocdepth}{1}
	
	\tableofcontents
	
	\section{Introduction} \label{sect_introduction}
	
	A convex body $K \subset \R^n$ is a convex compact set. Throughout this article, we assume (unless mentioned otherwise) that the convex bodies under consideration have non-empty interior. The \textit{covariance matrix} of a convex body $K \subset \R^n$ is given by
	\begin{equation}
		A_K \coloneqq \E[XX^\transp]-\E[X] \E[X]^\transp,
	\end{equation}
	where $X$ is uniformly distributed on $K$. The \textit{isotropic constant} $L_K$ of $K$ is then given by
	\begin{equation}
		L_K^{2n} = \frac{\det A_K}{(\vol K)^2}.
	\end{equation}
	Since $L_K$ is an affine invariant, we can restrict our attention to the case where $K$ is \textit{isotropic}, i.e., we have $\E[X] = 0$ and $A_K = I_n$.

	Until recently, it was a major open problem in high-dimensional convex geometry whether $L_K$ is bounded from above by a universal constant; the recently proven \textit{isotropic constant conjecture} asserts that such a constant exists. The problem was first considered by Bourgain \cite{bou86}. We refer to \cite{km22} for an overview of the theory of isotropic constants; a very detailed account is given in \cite{bgvv16}. Bourgain \cite{bou91} showed that $L_K \in O(\log n \cdot \sqrt[4]{n})$, which was improved by Klartag \cite{kla06} to $O(\sqrt[4]{n})$. More recently, Eldan's method of stochastic localization \cite{eld13} has led to a sequence of breakthroughs. In a seminal paper \cite{che21}, Chen showed that $L_K \in O(n^\varepsilon)$ for every $\varepsilon>0$, which was improved by Klartag and Lehec \cite{kl22} to the polylogarithmic bound $L_K \in O(\log^4 n)$, and then by Klartag \cite{kla23} to $L_K \in O(\sqrt{\log n})$. In December 2024, the isotropic constant conjecture was finally proven by Klartag and Lehec~\cite{kl24}, building on a preprint by Guan \cite{gua24}, which had been published on \texttt{arxiv.org} only a few days before.
	
	With regard to the natural follow-up question for a sharp upper bound on $L_K$ in any given dimension $n \geq 1$, a strong version of the isotropic constant conjecture asserts that every convex body $K \subset \R^n$ satisfies
	\begin{equation} \label{eq_strong_isotropic_constant_conjecture}
		L_K \leq L_{\Delta_n} = \frac{\sqrt[n]{n!}}{(n+1)^\frac{n+1}{2n} \cdot \sqrt{n+2}},
	\end{equation}
	where $\Delta_n$ is an $n$-dimensional simplex. %
	It was shown by Klartag \cite{kla18} that the strong isotropic constant conjecture implies the asymmetric version of the well-known Mahler conjecture.
	
	A convex body $K \subset \R^n$ is called \textit{centrally symmetric} if it has a center of symmetry, i.e., there exists a point $c \in \R^n$ with $K-c = -(K-c)$. A symmetric counterpart to \eqref{eq_strong_isotropic_constant_conjecture} is the conjecture that every centrally symmetric convex body $K \subset \R^n$ satisfies
	\begin{equation} \label{eq_symmetric_strong_isotropic_constant_conjecture}
		L_K \leq L_{C_n} = \frac{1}{\sqrt{12}},
	\end{equation}
	where $C_n$ is an $n$-dimensional parallelepiped.%
	
	A \textit{polytope} $P \subset\R^n$ is the convex hull of finitely many points. Throughout this article, we restrict our attention to polytopes that have non-empty interior.  A polytope is called \textit{simplicial} if every facet of $P$ is a simplex. We say that a vertex $v \in P$ is \textit{simplicial} if every facet $F \subset P$ with $v \in F$ is a simplex.
	
	In \cite{rad16}, Rademacher proved the following remarkable result, which provides supporting evidence for \eqref{eq_strong_isotropic_constant_conjecture}.
	
	\begin{theorem}[Rademacher] \label{thm_rademacher}
		Let $P \subset \R^n$ be a polytope that locally maximizes the isotropic constant. If $P$ is simplicial, then $P$ is a simplex.
	\end{theorem}
	
	The proof strategy of Theorem \ref{thm_rademacher} presented in \cite{rad16} also forms the basis of the present article. Broadly speaking, Theorem \ref{thm_rademacher} can be interpreted as saying that a polytopal local maximizer of $K \mapsto L_K$ with ``generic'' boundary structure has to be highly regular.  In the following, we present several new results in this spirit. As we vary the exact assumptions on $P$, each of the three high-dimensional ``Platonic solids'' (regular $n$-polytopes) will make its appearance.%
	
	We say that two facets $F_1, F_2$ of an $n$-dimensional polytope $P$ are \textit{adjacent} if the face $F_1 \cap F_2$ is a \textit{ridge} of $P$, i.e., an $(n-2)$-dimensional face of $P$. We note that a simplicial polytope $P$ has the following property: for every pair of adjacent facets $F_1,F_2$, there exists an affine map that maps $F_1$ to $F_2$ and fixes $F_1 \cap F_2$ pointwise. In Section \ref{sect_first_order_conditions}, we show that a polytopal maximizer of the isotropic constant with this property must be a simplex, generalizing Theorem \ref{thm_rademacher}.
	
	The proof of Theorem \ref{thm_rademacher} given in \cite{rad16} consists of two parts: using a first-order condition for the facets, it is first shown that the local maximizer $P$ must be symmetric with respect to a hyperplane $H$; then a result due to Campi, Colesanti and Gronchi \cite[Thm.~3.6]{ccg99} is used to show that it cannot be a local maximizer unless it is a simplex. Using the techniques from \cite{ccg99}, we show in Section \ref{sect_simplicial_vertices} that a polytope $P$ with a hyperplane symmetry ``locally at a simplicial vertex'' cannot be a maximizer, leading to the following alternative strengthening of Theorem \ref{thm_rademacher}. 
	
	\begin{theorem} \label{thm_simplicial_vertex}
		Let $P \subset \R^n$ be a polytope that locally maximizes the isotropic constant. If $P$ has a simplicial vertex, then $P$ is a simplex.
	\end{theorem}
	With regard to \eqref{eq_symmetric_strong_isotropic_constant_conjecture}, we show the following symmetric counterpart of Theorem \ref{thm_simplicial_vertex} in Section \ref{sect_the_centrally_symmetric_case}.
	\begin{theorem} \label{thm_symmetric_simplicial_vertex}
		Let $P \subset \R^n$ be a centrally symmetric polytope that locally maximizes the isotropic constant in the class of centrally symmetric convex bodies. If $P$ has a simplicial vertex, then $P$ is a cross-polytope.
	\end{theorem}
	
	Since an $n$-dimensional cross-polytope $C_n^*$ satisfies $L_{C_n^*}=\frac{\sqrt[n]{n!}}{\sqrt{2}\cdot \sqrt{n+1} \cdot \sqrt{n+2}} \leq \frac{1}{\sqrt{12}}$, this finding can be interpreted as supporting evidence for \eqref{eq_symmetric_strong_isotropic_constant_conjecture}.
	
	\begin{figure}[ht]
		\begin{subfigure}[b]{.33\linewidth}
			\centering
			\begin{tikzpicture}%
	[x={(0.562037cm, -0.330669cm)},
	y={(0.827112cm, 0.224669cm)},
	z={(0.000024cm, 0.916614cm)},
	scale=1.200000,
	back/.style={dotted, thin},
	edge/.style={color=black},
	facet/.style={fill=white,fill opacity=0.000000},
	vertex/.style={},
	back2/.style={dotted},
	edge2/.style={color=red,dashed},
	facet2/.style={fill=Emerald,fill opacity=0.20000},
	vertex2/.style={inner sep=1pt,circle,draw=darkgray!25!black,fill=darkgray!75!black,thick},
	back3/.style={dotted, thin},
	edge3/.style={color=blue,very thick},
	facet3/.style={fill=green,fill opacity=0.20000},
	vertex3/.style={inner sep=1pt,circle,draw=darkgray!25!black},
	vertex4/.style={inner sep=1pt,circle,draw=blue,fill=blue},
	facet4/.style={fill=blue,fill opacity=0.20000},]
	
	\coordinate (1.00000, -1.05000, 1.00000) at (1.00000, -1.05000, 1.00000);
	\coordinate (1.00000, 1.00000, -1.00000) at (1.00000, 1.00000, -1.00000);
	\coordinate (1.00000, 1.00000, 1.00000) at (1.00000, 1.00000, 1.00000);
	\coordinate (-1.00000, 1.00000, 1.00000) at (-1.00000, 1.00000, 1.00000);
	\coordinate (-1.00000, 1.00000, -1.00000) at (-1.00000, 1.00000, -1.00000);
	\coordinate (1.00000, -0.95000, -1.00000) at (1.00000, -0.95000, -1.00000);
	\coordinate (-1.00000, -1.05000, -1.00000) at (-1.00000, -1.05000, -1.00000);
	\coordinate (-1.00000, -0.95000, 1.00000) at (-1.00000, -0.95000, 1.00000);
	\coordinate (0.00000, -1.30000, 0.00000) at (0.00000, -1.30000, 0.00000);
	\draw[edge,back] (1.00000, 1.00000, -1.00000) -- (-1.00000, 1.00000, -1.00000);
	\draw[edge,back] (-1.00000, 1.00000, 1.00000) -- (-1.00000, 1.00000, -1.00000);
	\draw[edge,back] (-1.00000, 1.00000, -1.00000) -- (-1.00000, -1.05000, -1.00000);
	\node[vertex] at (-1.00000, 1.00000, -1.00000)     {};
	\fill[facet2] (0.00000, -1.30000, 0.00000) -- (1.00000, -0.95000, -1.00000) -- (-1.00000, -1.05000, -1.00000) -- cycle {};
	\fill[facet2] (0.00000, -1.30000, 0.00000) -- (-1.00000, -1.05000, -1.00000) -- (-1.00000, -0.95000, 1.00000) -- cycle {};
	\fill[facet2] (0.00000, -1.30000, 0.00000) -- (1.00000, -1.05000, 1.00000) -- (-1.00000, -0.95000, 1.00000) -- cycle {};
	\fill[facet2] (0.00000, -1.30000, 0.00000) -- (1.00000, -1.05000, 1.00000) -- (1.00000, -0.95000, -1.00000) -- cycle {};
	\draw[edge] (1.00000, -1.05000, 1.00000) -- (1.00000, 1.00000, 1.00000);
	\draw[edge] (1.00000, -1.05000, 1.00000) -- (1.00000, -0.95000, -1.00000);
	\draw[edge] (1.00000, -1.05000, 1.00000) -- (-1.00000, -0.95000, 1.00000);
	\draw[edge] (1.00000, -1.05000, 1.00000) -- (0.00000, -1.30000, 0.00000);
	\draw[edge] (1.00000, 1.00000, -1.00000) -- (1.00000, 1.00000, 1.00000);
	\draw[edge] (1.00000, 1.00000, -1.00000) -- (1.00000, -0.95000, -1.00000);
	\draw[edge] (1.00000, 1.00000, 1.00000) -- (-1.00000, 1.00000, 1.00000);
	\draw[edge] (-1.00000, 1.00000, 1.00000) -- (-1.00000, -0.95000, 1.00000);
	\draw[edge] (1.00000, -0.95000, -1.00000) -- (-1.00000, -1.05000, -1.00000);
	\draw[edge] (1.00000, -0.95000, -1.00000) -- (0.00000, -1.30000, 0.00000);
	\draw[edge] (-1.00000, -1.05000, -1.00000) -- (-1.00000, -0.95000, 1.00000);
	\draw[edge] (-1.00000, -1.05000, -1.00000) -- (0.00000, -1.30000, 0.00000);
	\draw[edge] (-1.00000, -0.95000, 1.00000) -- (0.00000, -1.30000, 0.00000);
	\node[vertex] at (1.00000, -1.05000, 1.00000)     {};
	\node[vertex] at (1.00000, 1.00000, -1.00000)     {};
	\node[vertex] at (1.00000, 1.00000, 1.00000)     {};
	\node[vertex] at (-1.00000, 1.00000, 1.00000)     {};
	\node[vertex] at (1.00000, -0.95000, -1.00000)     {};
	\node[vertex] at (-1.00000, -1.05000, -1.00000)     {};
	\node[vertex] at (-1.00000, -0.95000, 1.00000)     {};
	\node[vertex2] at (0.00000, -1.30000, 0.00000)     {};
\end{tikzpicture}
			
			\caption{simplicial vertex}
		\end{subfigure}
		\begin{subfigure}[b]{.32\linewidth}
			\centering
			\begin{tikzpicture}%
	[x={(0.562037cm, -0.330669cm)},
	y={(0.827112cm, 0.224669cm)},
	z={(0.000024cm, 0.916614cm)},
	scale=1.200000,
	back/.style={dotted, thin},
	edge/.style={color=black},
	facet/.style={fill=white,fill opacity=0.000000},
	vertex/.style={},
	back2/.style={dotted},
	edge2/.style={color=red,dashed},
	facet2/.style={fill=Emerald,fill opacity=0.20000},
	vertex2/.style={inner sep=1pt,circle,draw=darkgray!25!black,fill=darkgray!75!black,thick},
	back3/.style={dotted, thin},
	edge3/.style={color=black,very thick},
	facet3/.style={fill=Emerald,fill opacity=0.10000},
	vertex3/.style={inner sep=1pt,circle,draw=darkgray!25!black},
	facet4/.style={fill=blue,fill opacity=0.20000},]
	
	\coordinate (1.00000, -1.05000, 1.00000) at (1.00000, -1.05000, 1.00000);
	\coordinate (-1.00000, -1.05000, -1.00000) at (-1.00000, -1.05000, -1.00000);
	\coordinate (-1.00000, -0.95000, 1.00000) at (-1.00000, -0.95000, 1.00000);
	\coordinate (1.00000, -0.95000, -1.00000) at (1.00000, -0.95000, -1.00000);
	\coordinate (0.00000, -1.30000, 0.00000) at (0.00000, -1.30000, 0.00000);
	\coordinate (-1.00000, 1.05000, -1.00000) at (-1.00000, 1.05000, -1.00000);
	\coordinate (1.00000, 1.05000, 1.00000) at (1.00000, 1.05000, 1.00000);
	\coordinate (1.00000, 0.95000, -1.00000) at (1.00000, 0.95000, -1.00000);
	\coordinate (-1.00000, 0.95000, 1.00000) at (-1.00000, 0.95000, 1.00000);
	\coordinate (0.00000, 1.30000, 0.00000) at (0.00000, 1.30000, 0.00000);
	\draw[edge,back] (-1.00000, -1.05000, -1.00000) -- (-1.00000, 1.05000, -1.00000);
	\draw[edge,back] (-1.00000, 1.05000, -1.00000) -- (1.00000, 0.95000, -1.00000);
	\draw[edge,back] (-1.00000, 1.05000, -1.00000) -- (-1.00000, 0.95000, 1.00000);
	\draw[edge,back] (-1.00000, 1.05000, -1.00000) -- (0.00000, 1.30000, 0.00000);
	\draw[edge,back] (1.00000, 1.05000, 1.00000) -- (0.00000, 1.30000, 0.00000);
	\draw[edge,back] (1.00000, 0.95000, -1.00000) -- (0.00000, 1.30000, 0.00000);
	\draw[edge,back] (-1.00000, 0.95000, 1.00000) -- (0.00000, 1.30000, 0.00000);
	\node[vertex] at (-1.00000, 1.05000, -1.00000)     {};
	\node[vertex2] at (0.00000, 1.30000, 0.00000)     {};
	
	\fill[facet3] (0.00000, 1.30000, 0.00000) -- (-1.00000, 1.05000, -1.00000) -- (1.00000, 0.95000, -1.00000) -- cycle {};
	\fill[facet3] (0.00000, 1.30000, 0.00000) -- (1.00000, 1.05000, 1.00000) -- (1.00000, 0.95000, -1.00000) -- cycle {};
	\fill[facet3] (0.00000, 1.30000, 0.00000) -- (1.00000, 1.05000, 1.00000) -- (-1.00000, 0.95000, 1.00000) -- cycle {};
	\fill[facet3] (0.00000, 1.30000, 0.00000) -- (-1.00000, 1.05000, -1.00000) -- (-1.00000, 0.95000, 1.00000) -- cycle {};
	
	\fill[facet2] (0.00000, -1.30000, 0.00000) -- (1.00000, -1.05000, 1.00000) -- (-1.00000, -0.95000, 1.00000) -- cycle {};
	\fill[facet2] (0.00000, -1.30000, 0.00000) -- (-1.00000, -1.05000, -1.00000) -- (-1.00000, -0.95000, 1.00000) -- cycle {};
	\fill[facet2] (0.00000, -1.30000, 0.00000) -- (-1.00000, -1.05000, -1.00000) -- (1.00000, -0.95000, -1.00000) -- cycle {};
	\fill[facet2] (0.00000, -1.30000, 0.00000) -- (1.00000, -1.05000, 1.00000) -- (1.00000, -0.95000, -1.00000) -- cycle {};
	\draw[edge] (1.00000, -1.05000, 1.00000) -- (-1.00000, -0.95000, 1.00000);
	\draw[edge] (1.00000, -1.05000, 1.00000) -- (1.00000, -0.95000, -1.00000);
	\draw[edge] (1.00000, -1.05000, 1.00000) -- (0.00000, -1.30000, 0.00000);
	\draw[edge] (1.00000, -1.05000, 1.00000) -- (1.00000, 1.05000, 1.00000);
	\draw[edge] (-1.00000, -1.05000, -1.00000) -- (-1.00000, -0.95000, 1.00000);
	\draw[edge] (-1.00000, -1.05000, -1.00000) -- (1.00000, -0.95000, -1.00000);
	\draw[edge] (-1.00000, -1.05000, -1.00000) -- (0.00000, -1.30000, 0.00000);
	\draw[edge] (-1.00000, -0.95000, 1.00000) -- (0.00000, -1.30000, 0.00000);
	\draw[edge] (-1.00000, -0.95000, 1.00000) -- (-1.00000, 0.95000, 1.00000);
	\draw[edge] (1.00000, -0.95000, -1.00000) -- (0.00000, -1.30000, 0.00000);
	\draw[edge] (1.00000, -0.95000, -1.00000) -- (1.00000, 0.95000, -1.00000);
	\draw[edge] (1.00000, 1.05000, 1.00000) -- (1.00000, 0.95000, -1.00000);
	\draw[edge] (1.00000, 1.05000, 1.00000) -- (-1.00000, 0.95000, 1.00000);
	\node[vertex] at (1.00000, -1.05000, 1.00000)     {};
	\node[vertex] at (-1.00000, -1.05000, -1.00000)     {};
	\node[vertex] at (-1.00000, -0.95000, 1.00000)     {};
	\node[vertex] at (1.00000, -0.95000, -1.00000)     {};
	\node[vertex2] at (0.00000, -1.30000, 0.00000)     {};
	\node[vertex] at (1.00000, 1.05000, 1.00000)     {};
	\node[vertex] at (1.00000, 0.95000, -1.00000)     {};
	\node[vertex] at (-1.00000, 0.95000, 1.00000)     {};
\end{tikzpicture}
			
			\caption{symmetric, simplical vertex}
		\end{subfigure}
		\begin{subfigure}[b]{.32\linewidth}
			\centering
			\begin{tikzpicture}%
	[x={(0.562037cm, -0.330669cm)},
	y={(0.827112cm, 0.224669cm)},
	z={(0.000024cm, 0.916614cm)},
	scale=0.6750000,
	back/.style={dotted, thin},
	edge/.style={color=black},
	facet/.style={fill=white,fill opacity=0.000000},
	vertex/.style={},
	back2/.style={dotted},
	edge2/.style={color=red,dashed},
	facet2/.style={fill=Emerald,fill opacity=0.20000},
	vertex2/.style={inner sep=1pt,circle,draw=darkgray!25!black,fill=darkgray!75!black,thick},
	back3/.style={dotted, thin},
	edge3/.style={color=black,very thick},
	facet3/.style={fill=Emerald,fill opacity=0.10000},
	vertex3/.style={inner sep=1pt,circle,draw=darkgray!25!black},
	facet4/.style={fill=blue,fill opacity=0.20000},]
	
	\coordinate (-2.00000, -1.00000, 1.00000) at (-2.00000, -1.00000, 1.00000);
	\coordinate (-1.00000, -2.00000, -2.00000) at (-1.00000, -2.00000, -2.00000);
	\coordinate (2.00000, -1.00000, 1.00000) at (2.00000, -1.00000, 1.00000);
	\coordinate (2.00000, -1.00000, -1.00000) at (2.00000, -1.00000, -1.00000);
	\coordinate (2.00000, 1.00000, -1.00000) at (2.00000, 1.00000, -1.00000);
	\coordinate (1.00000, 2.00000, 2.00000) at (1.00000, 2.00000, 2.00000);
	\coordinate (1.00000, 2.00000, 0.00000) at (1.00000, 2.00000, 0.00000);
	\coordinate (2.00000, 1.00000, 1.00000) at (2.00000, 1.00000, 1.00000);
	\coordinate (1.00000, 0.00000, -2.00000) at (1.00000, 0.00000, -2.00000);
	\coordinate (1.00000, 0.00000, 2.00000) at (1.00000, 0.00000, 2.00000);
	\coordinate (1.00000, -2.00000, -2.00000) at (1.00000, -2.00000, -2.00000);
	\coordinate (-1.00000, 2.00000, 0.00000) at (-1.00000, 2.00000, 0.00000);
	\coordinate (1.00000, -2.00000, 0.00000) at (1.00000, -2.00000, 0.00000);
	\coordinate (-1.00000, 0.00000, -2.00000) at (-1.00000, 0.00000, -2.00000);
	\coordinate (-1.00000, 0.00000, 2.00000) at (-1.00000, 0.00000, 2.00000);
	\coordinate (-2.00000, 1.00000, 1.00000) at (-2.00000, 1.00000, 1.00000);
	\coordinate (-2.00000, 1.00000, -1.00000) at (-2.00000, 1.00000, -1.00000);
	\coordinate (-1.00000, 2.00000, 2.00000) at (-1.00000, 2.00000, 2.00000);
	\coordinate (-2.00000, -1.00000, -1.00000) at (-2.00000, -1.00000, -1.00000);
	\coordinate (-1.00000, -2.00000, 0.00000) at (-1.00000, -2.00000, 0.00000);
	\draw[edge,back] (-2.00000, -1.00000, 1.00000) -- (-2.00000, 1.00000, 1.00000);
	\draw[edge,back] (-2.00000, -1.00000, 1.00000) -- (-2.00000, -1.00000, -1.00000);
	\draw[edge,back] (-1.00000, -2.00000, -2.00000) -- (-1.00000, 0.00000, -2.00000);
	\draw[edge,back] (-1.00000, -2.00000, -2.00000) -- (-2.00000, -1.00000, -1.00000);
	\draw[edge,back] (1.00000, 2.00000, 0.00000) -- (-1.00000, 2.00000, 0.00000);
	\draw[edge,back] (1.00000, 0.00000, -2.00000) -- (-1.00000, 0.00000, -2.00000);
	\draw[edge,back] (-1.00000, 2.00000, 0.00000) -- (-2.00000, 1.00000, -1.00000);
	\draw[edge,back] (-1.00000, 2.00000, 0.00000) -- (-1.00000, 2.00000, 2.00000);
	\draw[edge,back] (-1.00000, 0.00000, -2.00000) -- (-2.00000, 1.00000, -1.00000);
	\draw[edge,back] (-2.00000, 1.00000, 1.00000) -- (-2.00000, 1.00000, -1.00000);
	\draw[edge,back] (-2.00000, 1.00000, 1.00000) -- (-1.00000, 2.00000, 2.00000);
	\draw[edge,back] (-2.00000, 1.00000, -1.00000) -- (-2.00000, -1.00000, -1.00000);
	\node[vertex] at (-1.00000, 2.00000, 0.00000)     {};
	\node[vertex] at (-1.00000, 0.00000, -2.00000)     {};
	\node[vertex] at (-2.00000, 1.00000, -1.00000)     {};
	\node[vertex] at (-2.00000, -1.00000, -1.00000)     {};
	\node[vertex] at (-2.00000, 1.00000, 1.00000)     {};
	\fill[facet3]  (-2.00000, -1.00000, -1.00000)  -- (-2.00000, -1.00000, 1.00000) -- (-1.00000, -2.00000, 0.00000) -- (-1.00000, -2.00000, -2.00000) -- cycle {};
	
	\fill[facet3]  (-2.00000, -1.00000, -1.00000)  -- (-2.00000, -1.00000, 1.00000) -- (-2.00000, 1.00000, 1.00000) -- (-2.00000, 1.00000, -1.00000) -- cycle {};
	
	\fill[facet3]  (-1.00000, 2.00000, 0.00000) -- (-1.00000, 2.00000, 2.00000)  -- (-2.00000, 1.00000, 1.00000) -- (-2.00000, 1.00000, -1.00000) -- cycle {};
	
	\fill[facet3] (1.00000, 2.00000, 2.00000) -- (-1.00000, 2.00000, 2.00000) -- (-1.00000, 2.00000, 0.00000) -- (1.00000, 2.00000, 0.00000)   -- cycle {};

	\fill[facet2] (1.00000, 2.00000, 2.00000) -- (2.00000, 1.00000, 1.00000) -- (2.00000, 1.00000, -1.00000) -- (1.00000, 2.00000, 0.00000) -- cycle {};
	\fill[facet2] (2.00000, 1.00000, 1.00000) -- (2.00000, -1.00000, 1.00000) -- (2.00000, -1.00000, -1.00000) -- (2.00000, 1.00000, -1.00000) -- cycle {};
	\fill[facet2] (1.00000, -2.00000, 0.00000) -- (2.00000, -1.00000, 1.00000) -- (2.00000, -1.00000, -1.00000) -- (1.00000, -2.00000, -2.00000) -- cycle {};
	\fill[facet2] (-1.00000, -2.00000, 0.00000) -- (-1.00000, -2.00000, -2.00000) -- (1.00000, -2.00000, -2.00000) -- (1.00000, -2.00000, 0.00000) -- cycle {};
	\draw[edge] (-2.00000, -1.00000, 1.00000) -- (-1.00000, 0.00000, 2.00000);
	\draw[edge] (-2.00000, -1.00000, 1.00000) -- (-1.00000, -2.00000, 0.00000);
	\draw[edge] (-1.00000, -2.00000, -2.00000) -- (1.00000, -2.00000, -2.00000);
	\draw[edge] (-1.00000, -2.00000, -2.00000) -- (-1.00000, -2.00000, 0.00000);
	\draw[edge] (2.00000, -1.00000, 1.00000) -- (2.00000, -1.00000, -1.00000);
	\draw[edge] (2.00000, -1.00000, 1.00000) -- (2.00000, 1.00000, 1.00000);
	\draw[edge] (2.00000, -1.00000, 1.00000) -- (1.00000, 0.00000, 2.00000);
	\draw[edge] (2.00000, -1.00000, 1.00000) -- (1.00000, -2.00000, 0.00000);
	\draw[edge] (2.00000, -1.00000, -1.00000) -- (2.00000, 1.00000, -1.00000);
	\draw[edge] (2.00000, -1.00000, -1.00000) -- (1.00000, -2.00000, -2.00000);
	\draw[edge] (2.00000, 1.00000, -1.00000) -- (1.00000, 2.00000, 0.00000);
	\draw[edge] (2.00000, 1.00000, -1.00000) -- (2.00000, 1.00000, 1.00000);
	\draw[edge] (2.00000, 1.00000, -1.00000) -- (1.00000, 0.00000, -2.00000);
	\draw[edge] (1.00000, 2.00000, 2.00000) -- (1.00000, 2.00000, 0.00000);
	\draw[edge] (1.00000, 2.00000, 2.00000) -- (2.00000, 1.00000, 1.00000);
	\draw[edge] (1.00000, 2.00000, 2.00000) -- (1.00000, 0.00000, 2.00000);
	\draw[edge] (1.00000, 2.00000, 2.00000) -- (-1.00000, 2.00000, 2.00000);
	\draw[edge] (1.00000, 0.00000, -2.00000) -- (1.00000, -2.00000, -2.00000);
	\draw[edge] (1.00000, 0.00000, 2.00000) -- (-1.00000, 0.00000, 2.00000);
	\draw[edge] (1.00000, -2.00000, -2.00000) -- (1.00000, -2.00000, 0.00000);
	\draw[edge] (1.00000, -2.00000, 0.00000) -- (-1.00000, -2.00000, 0.00000);
	\draw[edge] (-1.00000, 0.00000, 2.00000) -- (-1.00000, 2.00000, 2.00000);
	\node[vertex] at (-2.00000, -1.00000, 1.00000)     {};
	\node[vertex] at (-1.00000, -2.00000, -2.00000)     {};
	\node[vertex] at (2.00000, -1.00000, 1.00000)     {};
	\node[vertex] at (2.00000, -1.00000, -1.00000)     {};
	\node[vertex] at (2.00000, 1.00000, -1.00000)     {};
	\node[vertex] at (1.00000, 2.00000, 2.00000)     {};
	\node[vertex] at (1.00000, 2.00000, 0.00000)     {};
	\node[vertex] at (2.00000, 1.00000, 1.00000)     {};
	\node[vertex] at (1.00000, 0.00000, -2.00000)     {};
	\node[vertex] at (1.00000, 0.00000, 2.00000)     {};
	\node[vertex] at (1.00000, -2.00000, -2.00000)     {};
	\node[vertex] at (1.00000, -2.00000, 0.00000)     {};
	\node[vertex] at (-1.00000, 0.00000, 2.00000)     {};
	\node[vertex] at (-1.00000, 2.00000, 2.00000)     {};
	\node[vertex] at (-1.00000, -2.00000, 0.00000)     {};
\end{tikzpicture}
			
			\caption{zonotope, cubical zone}
		\end{subfigure}
		\caption{The settings of Theorem \ref{thm_simplicial_vertex}, Theorem \ref{thm_symmetric_simplicial_vertex} and Theorem \ref{thm_cubical_zone}, respectively. We show that local maximizers of $K \mapsto L_K$ with these properties have to be a simplex, a cross-polytope and a cube, respectively.}
	\end{figure}
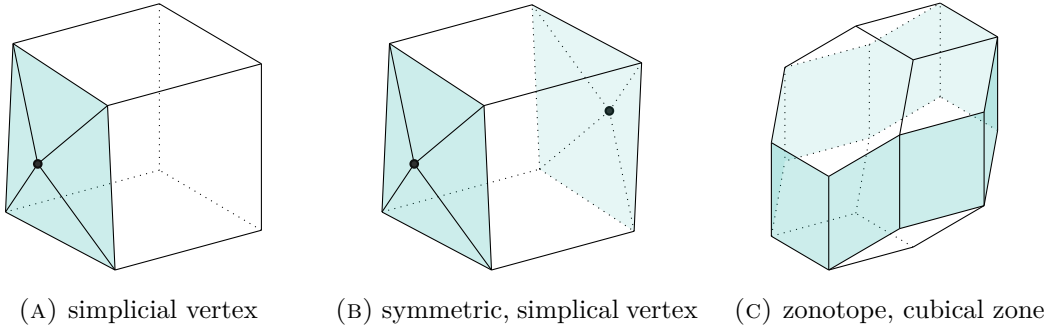
	
	Finally, we consider zonotopes. A \textit{zonotope} $Z$ is an affine image of a cube or, equivalently, the Minkowski sum of finitely many line segments, which are called the \textit{generators} of $Z$. A \textit{zonoid} is a convex body that can be approximated by zonotopes (with respect to the Hausdorff metric). The isotropic constant conjecture is known to be true for zonoids \cite{mp89}. It seems that the best known upper bound is $L_Z \leq \frac12$ \cite{pao03, bgvv16}. We note that the class of zonoids not only contains the conjectured maximizers from \eqref{eq_symmetric_strong_isotropic_constant_conjecture}, but also the minimizers in every dimension $n\geq 1$, which are known to be ellipsoids. 
	
	The set of all faces of a zonotope $Z$ that have a given generator as a Minkowski summand is called a \textit{zone} of $Z$. We say that a zone is \textit{cubical} if it contains only parallelepipeds. In Section \ref{sect_cubical_zones_of_zonotopes}, we show the following result, which can be interpreted as further supporting evidence for \eqref{eq_symmetric_strong_isotropic_constant_conjecture}.
	
	\begin{theorem} \label{thm_cubical_zone}
		Let $Z \subset \R^n$ be an isotropic zonotope that locally maximizes the isotropic constant in the class of centrally symmetric convex bodies. If $Z$ has a cubical zone, then $Z$ is a cube.
	\end{theorem}
	
	It was shown by McMullen \cite{mcm70} that if all ridges of a polytope $P$ are centrally symmetric, then $P$ is a zonotope; in particular, this holds if all facets of $P$ are parallelepipeds. Therefore, Theorem \ref{thm_cubical_zone} implies that a polytopal maximizer of $K \mapsto L_K$ in the class of centrally symmetric convex bodies must be a cube if its facets are all parallelepipeds.

	We emphasize that all theorems up to this point are hypothetical in the sense that it is unknown whether polytopes with the respective properties actually exist. In particular, Theorem \ref{thm_cubical_zone} yields no upper bound on isotropic constants of zonotopes. In Section \ref{sect_zonotopes_n+1_generators}, we turn to zonotopes with at most $n+1$ generators and show that \eqref{eq_symmetric_strong_isotropic_constant_conjecture} holds in this class. We also determine the minimizers in this class, which in three dimensions turn out to be affine images of the regular rhombic dodecahedron (and an $n$-dimensional analogue thereof in the general case).

	\section{First-order conditions for extremals of the isotropic constant} \label{sect_first_order_conditions}
	
	As mentioned above, Rademacher's proof of Theorem \ref{thm_rademacher} from \cite{rad16} makes use of two ingredients. The first ingredient is a first-order condition for the facets of a polytopal maximizer of the isotropic constant. Let $P \subset \R^n$ be an isotropic polytope. Intuitively speaking, if $P$ is a local maximizer of $K \mapsto L_K$, then the isotropic constant should not increase if some ``infinitesimal mass'' is added to or removed from $P$. It was shown in \cite{rad16} that if a facet $F \subset P$ is ``hinged'' around one of its ridges, then the derivative of $K \mapsto L_K^{2n}$ at $K=P$ is proportional to
	\begin{equation}
		\E[\norm{X}^2] - n-2,
	\end{equation}
	where the random vector $X$ is supported on $F$ and its density is proportional to the added ``infinitesimal mass'' (see also \cite{rad12} for a similar argument).
	
	In \cite{kip24}, the scope of such first-order arguments was studied in a somewhat more abstract framework. It was shown that if $P$ is an isotropic local maximizer of the isotropic constant, $F \subset P$ is a facet and $f \colon \relint F \rightarrow \R$ is an integrable concave function, then
	\begin{equation} \label{eq_spop_concave_condition}
		\int_F (\norm{x}_2^2-n-2) \cdot f(x) \dint \volt{n-1}(x) \leq 0
	\end{equation}
	holds \cite[Sect.~7]{kip24} and moreover, in some sense, all relevant first-order conditions on $P$ are implied by the conditions of type \eqref{eq_spop_concave_condition} \cite[Thm.~1.7]{kip24}. In particular, since the functions $x \mapsto x_i$ and $x \mapsto -x_i$ are concave for every $i \in [n]$, we get
	\begin{equation} \label{eq_spop_linear_condition}
		\int_F (\norm{x}_2^2-n-2) \cdot x_i \dint x = 0 \quad \text{for all } i \in [n].
	\end{equation}
	If $X$ is a random vector that is uniformly distributed on $F$, then \eqref{eq_spop_linear_condition} can be written more succinctly as
	\begin{equation} \label{eq_spop_condition}
		\E[\norm{X}^2 X] = (n+2) \E[X].
	\end{equation}
	In the remainder of this section, we will use \eqref{eq_spop_condition} to generalize Theorem \ref{thm_rademacher}.
	
	\begin{remark} \label{rem_there_are_no_polytopal_local_minimizers}
		If an isotropic polytope $P \subset \R^n$ is assumed to be a local minimizer of $K \mapsto L_K$, then the conditions of type \eqref{eq_spop_concave_condition} have to hold with the inequality sign reversed. In fact, this implies that such a polytope cannot exist. To see this, let $F \subset P$ be a facet and let $X$ be uniformly distributed on $F$. Now \eqref{eq_spop_condition} implies $\E[\norm{X}^2] = n+2$. It follows that $F_+ \coloneqq \{x \in F \mid \norm{x}> n+2\}$ satisfies $\volt{n-1}F_+>0$. Clearly, the function $f \colon \relint F \rightarrow \R$, $x \mapsto - \max\{\norm{x}^2_2,0\}$ is concave, but we have
		\begin{equation}
			\int_F (\norm{x}_2^2-n-2) \cdot f(x) \dint \volt{n-1}(x) =\int_{F_+} (\norm{x}_2^2-n-2) \cdot (-\norm{x}^2_2) \dint \volt{n-1}(x) <0,
		\end{equation}
		in contradiction to the reversed form of \eqref{eq_spop_concave_condition}.
	\end{remark}

	Let $K \subset \R^n$ be an $(n-1)$-dimensional convex body. A \textit{pyramid over $K$} is a convex body of the form $\conv( K \cup \{p\})$ for some point $p \notin \aff K$. For example, a simplex is a pyramid over each of its facets. In fact, the converse also holds. To see this, let $P$ be a polytope that is a pyramid over each of its facets. Because this property is preserved if we pass from $P$ to one of its facets, it follows by induction over the dimension that $P$ is a simplex. We will make use of this fact below.
	
	Let $P \subset \R^n$ be a polytope whose centroid is at the origin and let $F_1,F_2 \subset P$ be two adjacent facets. If $F_1$ and $F_2$ are pyramids over $F_1 \cap F_2$, then there exists a linear map $f$ that fixes $F_1 \cap F_2$ pointwise and maps $F_1$ to $F_2$. This leads us to the following definition, which was already alluded to in the introduction.
	
	\begin{definition} \label{def_reflector}
		Let $P \subset \R^n$ be a polytope, $G \subset P$ a ridge and $F_1,F_2 \subset P$ the two facets that contain $G$. Let $c_P$ denote the centroid of $P$. The ridge $G \subset P$ is called an \textit{affine reflector} if there exists an affine map $f\colon \R^n \rightarrow \R^n$ with
		\begin{equation}
			f(c_P)=c_P, \quad f(F_1) = F_2 \quad \text{and} \quad f|_G=\textstyle\id_G.
		\end{equation}
		An affine reflector $G$ is called a \textit{Euclidean reflector} if the corresponding map %
		$f$ is the reflection $\rho_G \colon \R^n \rightarrow \R^n$ across the hyperplane spanned by $G$ and the centroid of $P$.
	\end{definition}
	
	Considering the example above Definition \ref{def_reflector} further, we see that, in particular, every ridge of a simplicial polytope $P$ is an affine reflector. Using the first-order condition that was mentioned above, it was shown in \cite{rad16} that if $P$ is an isotropic simplicial polytope that is an extremal body of $K \mapsto L_K$, then every ridge of $P$ is even a Euclidean reflector. The following lemma generalizes this result to arbitrary affine reflectors.
	
	\begin{lemma} \label{lemma_euclidean_reflector}
		Let $P \subset \R^n$ be an isotropic polytope and let $G \subset P$ be an affine reflector. If every facet of $P$ satisfies \eqref{eq_spop_condition}, then $G$ is a Euclidean reflector.
	\end{lemma}
	
	\begin{proof}
		Let $F_1$, $F_2$ and $f$ be as in Definition \ref{def_reflector}. Since the centroid of $P$ is at the origin, $f$ is linear. Let $A\in \R^{n \times n}$ with $Ax=f(x)$ for all $x\in \R^n$. We note that $A$ is invertible because $\R^n=\linh F_2 \subset \im A$. %
		Let $X$ be the random vector that is uniformly distributed on $F_1$. Denoting the $i$-th standard basis vector of $\R^n$ by $e_i$, we assume without loss of generality that $\linh G = e_n^\bot$ and $F_1\subset \R^{n-1} \times [0,\infty)\subset \R^n$, which implies $\Prob(X_n>0) =1$. Because $f(F_1)=F_2$, the matrix $A$ is of the form $[e_1,\dots,e_{n-1},v]$ for some $v \in \R^n$. We have to show that $v=-e_n$. 
		
		Since $AX$ is uniformly distributed on $F_2$, \eqref{eq_spop_condition} implies
		\begin{equation}
			A\E[\norm{AX}^2 X]=\E[\norm{AX}^2 AX] = (n+2) \E[AX]= A(n+2) \E[X].
		\end{equation}
		Because $A$ is invertible, this is equivalent to
		\begin{equation}
			\E[\norm{AX}^2 X_i] = (n+2) \E[X_i] \quad \text{for } i \in [n].
		\end{equation}
		Let $i \in [n]$. Using $A=[e_1,\dots,e_{n-1},v]$, we get
		\begin{equation}
			(n+2) \E[X_i] = \E[\norm{AX}^2 X_i] = \E\left[\sum_{j=1}^{n-1}(X_j+v_jX_n)^2X_i+(v_nX_n)^2X_i\right],%
		\end{equation}
		which is equivalent to
		\begin{equation} \label{eq_inhom_lin}
			(n+2) \E[X_i] -\sum_{j=1}^{n-1}\E[X_i X_j^2] = 2 \cdot \sum_{j=1}^{n-1} \E[X_i X_jX_n] \cdot v_j+ \E[X_iX_n^2] \cdot \norm{v}^2.
		\end{equation}
		Let $\tilde{v}\coloneqq(v_1,\dots,v_{n-1}, \frac12 \norm{v}^2)$ and  define $M \in \R^{n\times n}$ by $m_{ij}=\E[X_i X_jX_n]$ for $i,j \in [n]$. If we denote the left-hand side of \eqref{eq_inhom_lin} by $b_i$, $i \in \N$, then the system of equations given by \eqref{eq_inhom_lin} can be rewritten as
		\begin{equation}
			b = M \tilde{v}.
		\end{equation}
		Because the facet $F_1$ satisfies \eqref{eq_spop_condition}, $\tilde{v}=\frac{1}{2} e_n$ is a solution of this system. It remains to show that $\tilde{v}$ is the only solution. We have		
		\begin{equation}
			y^\transp M  y= \E[y^\transp X X_n X^\transp y]= \E[ X_n \cdot |X^\transp y|^2]  \quad \text{for } y \in \R^n.
		\end{equation}
		Recalling that $\Prob(X_n>0) =1$ and observing that $\dim (F_1 \cap y^\bot)\leq n-2$ implies $\Prob(X^\transp y = 0) =0$, we conclude that $y^\transp M  y>0$ for all $y \in \R^n \setminus\{0\}$. This shows that $M$ is invertible. It follows that $\tilde{v}$ is the only solution of \eqref{eq_inhom_lin}. Since $A \neq I_n$, we have $v_n = -e_n$.
	\end{proof}
	
	Evidently, the existence of Euclidean reflectors in the boundary of $P$ imposes strong restrictions on the structure of $P$. More specifically, if every pair of facets is connected by a path of Euclidean reflectors (i.e., a path in the graph of $P^*$, the polar of $P$), then $P$ is \textit{monohedral}, i.e., all facets of $P$ are congruent. Moreover, if every ridge of $P$ is a Euclidean reflector, then $P$ is even \textit{isohedral}, i.e., the (full) symmetry group of $P$ acts transitively on the facets of $P$. This follows from the following lemma, whose statement is used in \cite{rad16}. We spell out the simple proof for later reference.
	
	\begin{lemma} \label{lemma_symmetric_wrt_hyperplane}
		Let $P \subset \R^n$ be a polytope. If every ridge of $P$ is a Euclidean reflector, then $P$ is symmetric with respect to the hyperplane spanned by any ridge $G \subset P$ and the centroid of $P$.
	\end{lemma}
	
	\begin{proof}
		Let $G \subset P$ be a fixed ridge and $F_{1},F_{-1}$ be the facets that contain $G$. Moreover, let $F\subset P$ be an arbitrary facet. Since the graph of $P^*$ is connected, there exists a sequence of facets $F_1,F_2,\dots, F_m =F$ where $F_i$ and $F_{i+1}$ are adjacent for every $i \in [m-1]$. Because $G_i \coloneqq F_i \cap F_{i+1}$ is a Euclidean reflector, it follows that $F_{i+1}= \rho_{G_i}(F_i)$  for every $i \in [m-1]$. Moreover, since $F_{-1} =\rho_G(F_1)$, the set $G_{-1}\coloneqq\rho_G(F_1 \cap F_2)$ is a ridge of $F_{-1}$.  Setting $F_{-2} \coloneqq\rho_{G_{-1}}(F_{-1}) =\rho_G (F_2)$ and proceeding iteratively, we obtain a sequence of facets $F_{-1},F_{-2},\dots, F_{-m}$ with $F_{-i}=\rho_G(F_i)$ for all $i \in [m]$. In particular, $\rho_G(F)=F_{-m}$ is a facet of $P$.
	\end{proof}
	
	Combining Lemma \ref{lemma_euclidean_reflector} and Lemma \ref{lemma_symmetric_wrt_hyperplane}, we obtain the following result, which generalizes \cite[Thm.~1]{rad16}. Here and in the following, we call a hyperplane \textit{linear} if it contains the origin.
	
	\begin{theorem} \label{thm_affine_reflector_isohedral}
			Let $P$ be an isotropic polytope at which the isotropic constant has a local extremum. If every ridge $G \subset P$ is an affine reflector, then $P$ is symmetric with respect to the linear hyperplane spanned by any ridge and, in particular, isohedral.
	\end{theorem}
	
	As discussed in Remark \ref{rem_there_are_no_polytopal_local_minimizers}, there are no polytopal local minimizers of the isotropic constant. Consequently, the statement about local minimizers in Theorem \ref{thm_affine_reflector_isohedral} is vacuously true if the term ``local extremum'' is interpreted as ``local extremum in the class of convex bodies''. However, it is easy to see that \eqref{eq_spop_linear_condition} and \eqref{eq_spop_condition} are necessary conditions for a polytopal local extremal body $P$ \textit{in the class of polytopes with at most $f_{n-1}(P)$ facets}, where $f_{n-1}(P)$ denotes the number of facets of $P$. Interpreting the term ``local extremum'' in this way, we strengthen Theorem \ref{thm_affine_reflector_isohedral} and obtain a non-vacuous statement for local minimizers.
	
	Up to this point, we have only used the conditions \eqref{eq_spop_linear_condition} and \eqref{eq_spop_condition}, which cannot distinguish between local maximizers and local minimizers of the isotropic constant. This brings us to the second ingredient of the proof of Theorem \ref{thm_rademacher} used in \cite{rad16}, namely, the following result from \cite{ccg99}. 
	
	\begin{theorem}[Campi-Colesanti-Gronchi] \label{thm_ccg_hyperplane_symmetry}
		Let $H \subset \R^n$ be a linear hyperplane and $K\subset \R^n$ a convex body. Moreover, let $L \coloneqq H^\bot$ be the orthogonal complement of $H$. %
		If $K$ is symmetric with respect to $H$, then $K$ cannot be a local maximizer of the isotropic constant unless the function $x \mapsto \volt{1}[K \cap (L+x)]$ is affine on $K$.
	\end{theorem}
	
	In the present section, we use Theorem \ref{thm_ccg_hyperplane_symmetry} as a ``black box''. The geometric mechanism behind this result will be discussed in the next section. Adding Theorem \ref{thm_ccg_hyperplane_symmetry} as a second ingredient to our previous results, we arrive at the main result of this section.	
	
	\begin{theorem} \label{thm_linear_reflector_simplex}
		Let $P$ be a polytope that is a local maximizer of the isotropic constant. If every ridge $G \subset P$ is an affine reflector, then $P$ is a simplex.
	\end{theorem}
	
	\begin{proof}
		Without loss of generality, we assume that $P$ is isotropic. Let $F \subset P$ be an arbitrary facet and let $F_1,\dots,F_\ell$ be the facets that are adjacent to $F$. We fix an index $i \in [\ell]$. By Theorem \ref{thm_affine_reflector_isohedral}, $P$ is isohedral and symmetric with respect to the linear hyperplane spanned by the ridge $F \cap F_i$. Setting $L_i \coloneqq [\linh (F \cap F_i)]^\bot$, we deduce from Theorem \ref{thm_ccg_hyperplane_symmetry} that the function $x \mapsto \volt{1}[P \cap (x+L_i)]$ is affine on $P$. This implies that for every $x \in P$, the extreme points of the line segment $P \cap (x+L_i)$ are contained in $F \cup F_i$. Now if $v\in P \cap (x+L_i)$ is a vertex of $P$, then $v$ is an extreme point of $P \cap (x+L_i)$. It follows that $\vertices P \subset F \cup F_i$ and hence $(\vertices P) \setminus F \subset F_i$. Because $i$ was arbitrary, we even have
		\begin{equation}
			\emptyset \neq (\vertices P) \setminus F \subset F_1 \cap \dots \cap F_\ell.
		\end{equation}
		We claim that $(\vertices P) \setminus F$ contains exactly one element. For this, we consider the polar polytope $P^*$ of $P$ and write $Q^\diamond\subset P^*$ for the polar face of any face $Q \subset P$. %
		Because $F_1,\dots,F_\ell$ are exactly the facets adjacent to $F$, the vertices $F^\diamond_1,\dots,F^\diamond_\ell \in P^*$ are exactly the neighbors of $F^\diamond$, which implies that $\dim\aff\{F^\diamond_1,\dots,F^\diamond_\ell\}\geq n-1$. Since $F_1 \cap \dots \cap F_\ell$ is non-empty, the vertices $F^\diamond_1,\dots,F^\diamond_\ell$ are contained in a common proper face $Q$ of $P^*$, which is a uniquely determined facet because $\dim Q \geq\dim\aff\{F^\diamond_1,\dots,F^\diamond_\ell\}\geq n-1$. With regard to $P$, this means that the set $(\vertices P) \setminus F \subset F_1 \cap \dots \cap F_\ell$ consists of the single vertex $Q^\diamond$, as claimed. %
		In other words, $P$ is a pyramid over $F$. Since $F$ was arbitrary, it follows that $P$ is a pyramid over each of its facets and hence a simplex.
	\end{proof}
	
	\section{Simplicial vertices} \label{sect_simplicial_vertices}
	
	The goal of this section is to prove Theorem \ref{thm_simplicial_vertex}. We start with a lemma which asserts that a polytopal extremal body of $L_K \mapsto K$ with a simplicial vertex $v$ has a hyperplane symmetry ``locally at $v$''.
	
	\begin{lemma} \label{lemma_local_symmetry}
		Let $P \subset \R^n$ be an isotropic polytope that satisfies \eqref{eq_spop_condition}. If $P$ has a simplicial vertex $v$, then the set
		\begin{equation}
			S_v\coloneqq\bigcup\{F \subset P \mid F \text{ is a facet of } P \text{ with } v\in F \}
		\end{equation}
		is invariant under the reflection $\rho_G$ for every ridge $G \subset P$ with $v \in G$, i.e., we have
		\begin{equation}
			\rho_G(S_v) = S_v.
		\end{equation}
	\end{lemma}
	In the proof of Lemma \ref{lemma_local_symmetry}, we use the fact that the set
	\begin{equation}
		\{F \subset P  \mid F \text{ is a face of } P \text{ with } v\in F\},
	\end{equation}
	partially ordered by set inclusion, is isomorphic to the face lattice of an $(n-1)$-dimensional polytope, which is called the \textit{vertex figure of $P$ at v} \cite[Prop.~2.4]{zie07}. In particular, since the graph of any polytope is connected, any two facets of $P$ that contain $v$ are connected by a path of adjacent facets that contain $v$.
	
	\begin{proof}[Proof of Lemma \ref{lemma_local_symmetry}]
		Let $F\subset P$ be a facet with $v \in P$ and let $G \subset P$ be a ridge with $v \in G$. Moreover, let $F_1,F_{-1}$ be the two facets of $P$ that contain $G$. Considering the vertex figure of $P$ at $v$, we see that there exists a sequence of facets $F_1,F_2,\dots, F_m =F$ that contain $v$ and where $F_i$ and $F_{i+1}$ are adjacent for every $i \in [m-1]$. By Lemma \ref{lemma_euclidean_reflector}, $F_i \cap F_{i+1}$ is a Euclidean reflector for every $i \in [m-1]$. Applying the same argument as in the proof of Lemma \ref{lemma_symmetric_wrt_hyperplane}, we obtain that $\rho_G(F)$ is a facet of $P$ that contains $v$.
	\end{proof}
	
	We now come to the geometric mechanism behind Theorem \ref{thm_ccg_hyperplane_symmetry}. It is given by certain ``movements'' of a convex body $K$, which can be seen as a continuous generalization of Steiner symmetrization and Blaschke shaking. The following definition is adapted from \cite[Def.~2.2]{ccg99}.
	
	\begin{definition} \label{def_rs_movement}
		Let $K \subset \R^n$ be a convex body and $v \in \R^n \setminus \{0\}$. We denote the orthogonal projection $\R^n \rightarrow \{v\}^\bot$ onto the orthogonal complement of $v$ by $\pi_{v^\bot}$. Let $\beta \colon \pi_{v^\bot}(K) \rightarrow \R$ be a function.
		If $0 \in I \subset \R$ is an interval such that the set
		\begin{equation} \label{eq_def_rs_movement}
			K_t \coloneqq \{x+t \beta(\pi_{v^\bot}(x))v \mid x \in K\}
		\end{equation}
		is convex for every $t \in I$, then the family of convex bodies $(K_t)_{t \in I}$ is called an \textit{RS-movement} of $K$. The function $\beta$ is called \textit{speed function of $(K_t)_{t \in I}$}.
	\end{definition}
	
	The relevance of RS-movements for our purposes is given by the fact that the isotropic constant of an RS-movement is a quasiconvex function of the time parameter $t \in I$. We will state this fact formally in Theorem \ref{thm_ccg} after some further preliminaries.
	
	We say that a function $f$ \textit{factors through} another function $g$ if there is a function $h$ with $f = h \circ g$. In the following, we sometimes call functions $\beta \colon \R^n \rightarrow \R$ which factor through $\pi_{v^\bot}$ \textit{speed functions}. Formally, the speed function in the sense of Definition \ref{def_rs_movement} is then the restriction of $\beta$ to $\pi_{v^\bot}(K)$. Clearly, if $K \subset \R^n$ is an arbitrary convex body, then every affine function $\beta \colon \R^n \rightarrow \R$ trivially induces a RS-movement $(K_t)_{t \in \R}$ of $K$. The interesting question is whether $K$ admits non-trivial RS-movements. This leads us to the following definition.
	
	\begin{definition} \label{def_rs_decomposable}
		Following \cite{ccg99}, we say that a convex body $K \subset \R^n$ is \textit{RS-indecomposable} if every RS-movement $(K_t)_{t \in [-\varepsilon,\varepsilon]}$ with $K=K_0$ and $\varepsilon>0$ has an affine speed function.%
	\end{definition}
	
	We now state the crucial property of RS-movements formally. The following theorem is a special case of \cite[Thm.~3.1]{ccg99}, corresponding to the exponent $r=2$. The latter result is stated in the setting of Sylvester's problem and addresses arbitrary exponents. For the connection between Sylvester's problem and isotropic constants, we refer to \cite{kin69}.
	
	\begin{theorem}[Campi-Colesanti-Gronchi] \label{thm_ccg}
		For any RS-movement $(K_t)_{t \in [a,b]}$, the function $[a,b] \rightarrow \R, t \mapsto L_{K_t}^{2n}$ is convex. If the speed function $\beta$ is not affine, then $t \mapsto L_{K_t}^{2n}$ is strictly convex. In particular, every local maximizer of the isotropic constant is RS-indecomposable.
	\end{theorem}

	Theorem \ref{thm_ccg_hyperplane_symmetry} now follows immediately from Theorem \ref{thm_ccg} by considering the speed function $x \mapsto -\volt{1}[K \cap (L+x)]$.
	
	We have already shown that an isotropic polytope $P$ that maximizes the isotropic constant has a hyperplane symmetry locally at each of its simplicial vertices. The key observation that leads to Theorem \ref{thm_symmetric_simplicial_vertex} is that such a local symmetry is sufficient to construct a non-trivial RS-movement, unless $P$ is a simplex. In some sense, this construction generalizes the solution of the planar case of \eqref{eq_strong_isotropic_constant_conjecture} by Campi, Colesanti and Gronchi \cite{ccg99}.
	
	We need the following lemma, which is taken from \cite[Ex.~2.6]{ccg99}.
	
	\begin{lemma} \label{lemma_pyramid_rs_decomposable}
		A pyramid $\tilde{K} \subset \R^n$ over an $(n-1)$-dimensional convex body $K \subset \R^n$ is RS-indecomposable if and only if $K$ is RS-indecomposable.
	\end{lemma}
	
	\begin{proof}[Proof of Theorem \ref{thm_simplicial_vertex}]
		Without loss of generality, we assume that $P$ is isotropic. Let $v \in P$ be a simplicial vertex and let $S_v$ be as in Lemma \ref{lemma_local_symmetry}. Moreover, let $G \subset P$ be any ridge with $v \in G$ and let $H \coloneqq \linh G$. Because $P$ is a local maximizer of the isotropic constant, it satisfies \eqref{eq_spop_condition} and hence Lemma \ref{lemma_local_symmetry} implies that $S_v$ is symmetric with respect to $H$. We define $\tilde{\beta}\colon \bd P \rightarrow \R$ by setting $\tilde{\beta}(v)=1$,  $\tilde{\beta}(w)=0$ for all $w \in \vertices  P\setminus \{v\}$ and interpolating affinely on each facet. This is well-defined because every facet that contains $v$ is a simplex; on the other facets the function $\tilde{\beta}$ vanishes. Since $S_v$ is symmetric with respect to the hyperplane $H$, the function  $\tilde{\beta}$ factors through $\pi_{H}$, i.e., there exists $\beta \colon \pi_{H}(P) \rightarrow \R$ with $\tilde{\beta} = \beta \circ \pi_{H}$. Let $u \in \sphere$ be a normal vector of $H$ and set
		\begin{equation}
			P_t \coloneqq \{x+t \beta(\pi_{H}(x))u \mid x \in P\} \quad \text{for } t \in \R.
		\end{equation}
		\begin{figure}[b]
			\begin{subfigure}[b]{.33\linewidth}
				\centering
				\begin{tikzpicture}%
	[x={(0.562037cm, -0.330669cm)},
	y={(0.827112cm, 0.224669cm)},
	z={(0.000024cm, 0.916614cm)},
	scale=1.200000,
	back/.style={dotted, thin},
	edge/.style={color=black},
	facet/.style={fill=white,fill opacity=0.000000},
	vertex/.style={},
	back2/.style={dotted},
	edge2/.style={color=red,dashed},
	facet2/.style={fill=orange,fill opacity=0.20000},
	vertex2/.style={inner sep=1pt,circle,draw=darkgray!25!black,fill=darkgray!75!black,thick},
	back3/.style={dotted, thin},
	edge3/.style={color=blue,very thick},
	facet3/.style={fill=green,fill opacity=0.20000},
	vertex3/.style={inner sep=1pt,circle,draw=darkgray!25!black},
	vertex4/.style={inner sep=1pt,circle,draw=blue,fill=blue},
	facet4/.style={fill=blue,fill opacity=0.20000},]
	\coordinate (1.00000, -1.05000, 1.00000) at (1.00000, -1.05000, 1.00000);
	\coordinate (1.00000, 1.00000, -1.00000) at (1.00000, 1.00000, -1.00000);
	\coordinate (1.00000, 1.00000, 1.00000) at (1.00000, 1.00000, 1.00000);
	\coordinate (-1.00000, 1.00000, 1.00000) at (-1.00000, 1.00000, 1.00000);
	\coordinate (-1.00000, 1.00000, -1.00000) at (-1.00000, 1.00000, -1.00000);
	\coordinate (1.00000, -0.95000, -1.00000) at (1.00000, -0.95000, -1.00000);
	\coordinate (-1.00000, -1.05000, -1.00000) at (-1.00000, -1.05000, -1.00000);
	\coordinate (-1.00000, -0.95000, 1.00000) at (-1.00000, -0.95000, 1.00000);
	\coordinate (0.00000, -1.30000, 0.00000) at (0.00000, -1.30000, 0.00000);
	\draw[edge,back] (1.00000, 1.00000, -1.00000) -- (-1.00000, 1.00000, -1.00000);
	\draw[edge,back] (-1.00000, 1.00000, 1.00000) -- (-1.00000, 1.00000, -1.00000);
	\draw[edge,back] (-1.00000, 1.00000, -1.00000) -- (-1.00000, -1.05000, -1.00000);
	\node[vertex] at (-1.00000, 1.00000, -1.00000)     {};
	\fill[facet4] (1.00000, 1.00000, -1.00000) -- (-1.00000, 1.00000, 1.00000) -- (-1.00000, -1.3000, 1.0000) -- (1.00000, -1.3000, -1.00000)  -- cycle {};
	\draw[edge] (1.00000, -1.05000, 1.00000) -- (1.00000, 1.00000, 1.00000);
	\draw[edge] (1.00000, -1.05000, 1.00000) -- (1.00000, -0.95000, -1.00000);
	\draw[edge] (1.00000, -1.05000, 1.00000) -- (-1.00000, -0.95000, 1.00000);
	\draw[edge] (1.00000, -1.05000, 1.00000) -- (0.00000, -1.30000, 0.00000);
	\draw[edge] (1.00000, 1.00000, -1.00000) -- (1.00000, 1.00000, 1.00000);
	\draw[edge] (1.00000, 1.00000, -1.00000) -- (1.00000, -0.95000, -1.00000);
	\draw[edge] (1.00000, 1.00000, 1.00000) -- (-1.00000, 1.00000, 1.00000);
	\draw[edge] (-1.00000, 1.00000, 1.00000) -- (-1.00000, -0.95000, 1.00000);
	\draw[edge] (1.00000, -0.95000, -1.00000) -- (-1.00000, -1.05000, -1.00000);
	\draw[edge3] (1.00000, -0.95000, -1.00000) -- (0.00000, -1.30000, 0.00000);
	\draw[edge] (-1.00000, -1.05000, -1.00000) -- (-1.00000, -0.95000, 1.00000);
	\draw[edge] (-1.00000, -1.05000, -1.00000) -- (0.00000, -1.30000, 0.00000);
	\draw[edge] (-1.00000, -0.95000, 1.00000) -- (0.00000, -1.30000, 0.00000);
	
	\draw[edge2] (1.00000, -1.30000, 1.00000) -- (-1.00000, -1.30000, -1.00000);
	\node[vertex] at (1.00000, -1.05000, 1.00000)     {};
	\node[vertex] at (1.00000, 1.00000, -1.00000)     {};
	\node[vertex] at (1.00000, 1.00000, 1.00000)     {};
	\node[vertex] at (-1.00000, 1.00000, 1.00000)     {};
	\node[vertex] at (1.00000, -0.95000, -1.00000)     {};
	\node[vertex] at (-1.00000, -1.05000, -1.00000)     {};
	\node[vertex] at (-1.00000, -0.95000, 1.00000)     {};
	\node[vertex2] at (0.00000, -1.30000, 0.00000)     {};
	\node[vertex3] at (1.00000, -1.30000, 1.00000)     {};
	\node[vertex3] at (-1.00000, -1.30000, -1.00000)     {};
	\node[vertex4] at (0.00000, 0.0000, 0.00000)     {};
\end{tikzpicture}
				
				\caption{initial polytope $P_0$}
			\end{subfigure}
			\begin{subfigure}[b]{.32\linewidth}
				\centering
				\begin{tikzpicture}%
	[x={(0.562037cm, -0.330669cm)},
	y={(0.827112cm, 0.224669cm)},
	z={(0.000024cm, 0.916614cm)},
	scale=1.200000,
	back/.style={dotted, thin},
	edge/.style={color=black},
	facet/.style={fill=white,fill opacity=0.000000},
	vertex/.style={},
	back2/.style={dotted},
	edge2/.style={color=red,dashed},
	facet2/.style={fill=Emerald,fill opacity=0.20000},
	vertex2/.style={inner sep=1pt,circle,draw=darkgray!25!black,fill=darkgray!75!black,thick},
	back3/.style={dotted, thin},
	edge3/.style={color=black,very thick},
	facet3/.style={fill=green,fill opacity=0.20000},
	vertex3/.style={inner sep=1pt,circle,draw=darkgray!25!black},
	facet4/.style={fill=blue,fill opacity=0.20000},]
	\coordinate (1.00000, -0.95000, -1.00000) at (1.00000, -0.95000, -1.00000);
	\coordinate (1.00000, 1.00000, -1.00000) at (1.00000, 1.00000, -1.00000);
	\coordinate (1.00000, 1.00000, 1.00000) at (1.00000, 1.00000, 1.00000);
	\coordinate (-1.00000, 1.00000, 1.00000) at (-1.00000, 1.00000, 1.00000);
	\coordinate (-1.00000, 1.00000, -1.00000) at (-1.00000, 1.00000, -1.00000);
	\coordinate (1.00000, -1.05000, 1.00000) at (1.00000, -1.05000, 1.00000);
	\coordinate (-1.00000, -0.95000, 1.00000) at (-1.00000, -0.95000, 1.00000);
	\coordinate (-1.00000, -1.05000, -1.00000) at (-1.00000, -1.05000, -1.00000);
	\coordinate (-0.40000, -1.30000, -0.40000) at (-0.40000, -1.30000, -0.40000);
	\draw[edge,back] (1.00000, 1.00000, -1.00000) -- (-1.00000, 1.00000, -1.00000);
	\draw[edge,back] (-1.00000, 1.00000, 1.00000) -- (-1.00000, 1.00000, -1.00000);
	\draw[edge,back] (-1.00000, 1.00000, -1.00000) -- (-1.00000, -1.05000, -1.00000);
	\node[vertex] at (-1.00000, 1.00000, -1.00000)     {};
	\fill[facet2] (-0.40000, -1.30000, -0.40000) -- (1.00000, -0.95000, -1.00000) -- (-1.00000, -1.05000, -1.00000) -- cycle {};
	\fill[facet2] (-0.40000, -1.30000, -0.40000) -- (-1.00000, -0.95000, 1.00000) -- (-1.00000, -1.05000, -1.00000) -- cycle {};
	\fill[facet2] (-0.40000, -1.30000, -0.40000) -- (1.00000, -1.05000, 1.00000) -- (-1.00000, -0.95000, 1.00000) -- cycle {};
	\fill[facet2] (-0.40000, -1.30000, -0.40000) -- (1.00000, -0.95000, -1.00000) -- (1.00000, -1.05000, 1.00000) -- cycle {};
	\draw[edge] (1.00000, -0.95000, -1.00000) -- (1.00000, 1.00000, -1.00000);
	\draw[edge] (1.00000, -0.95000, -1.00000) -- (1.00000, -1.05000, 1.00000);
	\draw[edge] (1.00000, -0.95000, -1.00000) -- (-1.00000, -1.05000, -1.00000);
	\draw[edge] (1.00000, -0.95000, -1.00000) -- (-0.40000, -1.30000, -0.40000);
	\draw[edge] (1.00000, 1.00000, -1.00000) -- (1.00000, 1.00000, 1.00000);
	\draw[edge] (1.00000, 1.00000, 1.00000) -- (-1.00000, 1.00000, 1.00000);
	\draw[edge] (1.00000, 1.00000, 1.00000) -- (1.00000, -1.05000, 1.00000);
	\draw[edge] (-1.00000, 1.00000, 1.00000) -- (-1.00000, -0.95000, 1.00000);
	\draw[edge] (1.00000, -1.05000, 1.00000) -- (-1.00000, -0.95000, 1.00000);
	\draw[edge] (1.00000, -1.05000, 1.00000) -- (-0.40000, -1.30000, -0.40000);
	\draw[edge] (-1.00000, -0.95000, 1.00000) -- (-1.00000, -1.05000, -1.00000);
	\draw[edge] (-1.00000, -0.95000, 1.00000) -- (-0.40000, -1.30000, -0.40000);
	\draw[edge] (-1.00000, -1.05000, -1.00000) -- (-0.40000, -1.30000, -0.40000);
	
	\draw[edge2] (1.00000, -1.30000, 1.00000) -- (-1.00000, -1.30000, -1.00000);
	\node[vertex] at (1.00000, -0.95000, -1.00000)     {};
	\node[vertex] at (1.00000, 1.00000, -1.00000)     {};
	\node[vertex] at (1.00000, 1.00000, 1.00000)     {};
	\node[vertex] at (-1.00000, 1.00000, 1.00000)     {};
	\node[vertex] at (1.00000, -1.05000, 1.00000)     {};
	\node[vertex] at (-1.00000, -0.95000, 1.00000)     {};
	\node[vertex] at (-1.00000, -1.05000, -1.00000)     {};
	\node[vertex3] at (1.00000, -1.30000, 1.00000)     {};
	\node[vertex3] at (-1.00000, -1.30000, -1.00000)     {};
	\node[vertex2] at (-0.40000, -1.30000, -0.40000)     {};
	
\end{tikzpicture}
				
				\caption{$v$ moves in direction $u$}
			\end{subfigure}
			\begin{subfigure}[b]{.32\linewidth}
				\centering
				\begin{tikzpicture}%
	[x={(0.562037cm, -0.330669cm)},
	y={(0.827112cm, 0.224669cm)},
	z={(0.000024cm, 0.916614cm)},
	scale=1.200000,
	back/.style={dotted, thin},
	edge/.style={color=black},
	facet/.style={fill=white,fill opacity=0.000000},
	vertex/.style={},
	back2/.style={dotted},
	edge2/.style={color=red,dashed},
	facet2/.style={fill=Emerald,fill opacity=0.20000},
	vertex2/.style={inner sep=1pt,circle,draw=darkgray!25!black,fill=darkgray!75!black,thick},
	back3/.style={dotted, thin},
	edge3/.style={color=black,very thick},
	facet3/.style={fill=green,fill opacity=0.20000},
	vertex3/.style={inner sep=1pt,circle,draw=darkgray!25!black},
	facet4/.style={fill=blue,fill opacity=0.20000},]
	\coordinate (1.00000, -0.95000, -1.00000) at (1.00000, -0.95000, -1.00000);
	\coordinate (1.00000, 1.00000, -1.00000) at (1.00000, 1.00000, -1.00000);
	\coordinate (1.00000, 1.00000, 1.00000) at (1.00000, 1.00000, 1.00000);
	\coordinate (-1.00000, 1.00000, 1.00000) at (-1.00000, 1.00000, 1.00000);
	\coordinate (-1.00000, 1.00000, -1.00000) at (-1.00000, 1.00000, -1.00000);
	\coordinate (-1.00000, -0.95000, 1.00000) at (-1.00000, -0.95000, 1.00000);
	\coordinate (1.00000, -1.05000, 1.00000) at (1.00000, -1.05000, 1.00000);
	\coordinate (-1.00000, -1.30000, -1.00000) at (-1.00000, -1.30000, -1.00000);
	\draw[edge,back] (1.00000, 1.00000, -1.00000) -- (-1.00000, 1.00000, -1.00000);
	\draw[edge,back] (-1.00000, 1.00000, 1.00000) -- (-1.00000, 1.00000, -1.00000);
	\draw[edge,back] (-1.00000, 1.00000, -1.00000) -- (-1.00000, -1.30000, -1.00000);
	\node[vertex] at (-1.00000, 1.00000, -1.00000)     {};
	\fill[facet2] (-1.00000, -1.30000, -1.00000) -- (-1.00000, -0.95000, 1.00000) -- (1.00000, -1.05000, 1.00000) -- cycle {};
	\fill[facet2] (-1.00000, -1.30000, -1.00000) -- (1.00000, -0.95000, -1.00000) -- (1.00000, -1.05000, 1.00000) -- cycle {};
	\draw[edge] (1.00000, -0.95000, -1.00000) -- (1.00000, 1.00000, -1.00000);
	\draw[edge] (1.00000, -0.95000, -1.00000) -- (1.00000, -1.05000, 1.00000);
	\draw[edge] (1.00000, -0.95000, -1.00000) -- (-1.00000, -1.30000, -1.00000);
	\draw[edge] (1.00000, 1.00000, -1.00000) -- (1.00000, 1.00000, 1.00000);
	\draw[edge] (1.00000, 1.00000, 1.00000) -- (-1.00000, 1.00000, 1.00000);
	\draw[edge] (1.00000, 1.00000, 1.00000) -- (1.00000, -1.05000, 1.00000);
	\draw[edge] (-1.00000, 1.00000, 1.00000) -- (-1.00000, -0.95000, 1.00000);
	\draw[edge] (-1.00000, -0.95000, 1.00000) -- (1.00000, -1.05000, 1.00000);
	\draw[edge] (-1.00000, -0.95000, 1.00000) -- (-1.00000, -1.30000, -1.00000);
	\draw[edge] (1.00000, -1.05000, 1.00000) -- (-1.00000, -1.30000, -1.00000);
	
	\draw[edge2] (1.00000, -1.30000, 1.00000) -- (-1.00000, -1.30000, -1.00000);
	\node[vertex] at (1.00000, -0.95000, -1.00000)     {};
	\node[vertex] at (1.00000, 1.00000, -1.00000)     {};
	\node[vertex] at (1.00000, 1.00000, 1.00000)     {};
	\node[vertex] at (-1.00000, 1.00000, 1.00000)     {};
	\node[vertex] at (-1.00000, -0.95000, 1.00000)     {};
	\node[vertex] at (1.00000, -1.05000, 1.00000)     {};
	\node[vertex3] at (1.00000, -1.30000, 1.00000)     {};
	\node[vertex2] at (-1.00000, -1.30000, -1.00000)     {};
\end{tikzpicture}
				
				\caption{end of movement $P_\varepsilon$}
			\end{subfigure}
			\caption{The RS-movement $(P_t)_{i \in [-\varepsilon,\varepsilon]}$. The vertex $v$ moves orthogonally to the linear hyperplane spanned by the ridge $G$ and the origin (marked in blue).}
		\end{figure}
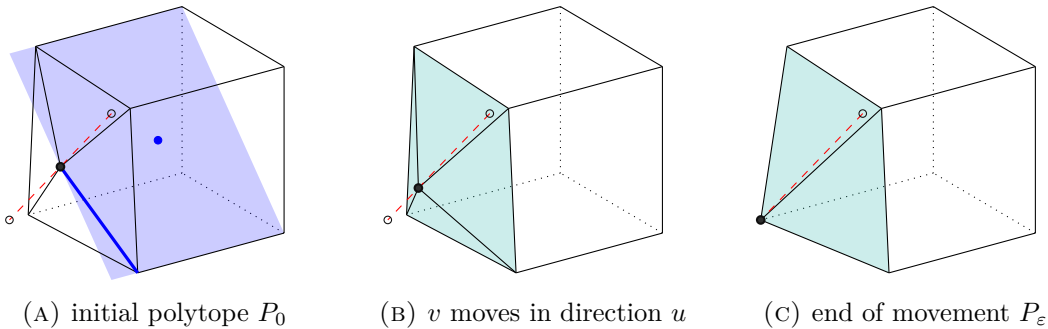%
		We set $\tilde{P} \coloneqq \conv \left(\vertices P \setminus \{v\}\right)$. Since $P$ has only finitely many vertices, there exists an $\varepsilon >0$ such that		
		\begin{equation}
			P_t = \conv ( \tilde{P} \cup \{v+tu\}) \quad \text{for all } t\in [-\varepsilon,\varepsilon].
		\end{equation}
		By Theorem \ref{thm_ccg}, it follows that $\beta$ is affine. This implies that $\pi_H(S_v) = \pi_H(P)$ and that there exists a hyperplane $H_v$ with $\vertices  P\setminus \{v\} \subset H_v$. In other words, $P$ is a pyramid over a polytope $Q \subset H_v$. Because every facet $F \subset P$ that contains $v$ is a simplex, $Q$ is simplicial. Let $w \in \vertices Q$. Since $G$ was arbitrary, we may assume without loss of generality that $w \in G$. Setting
		\begin{equation}
			\tilde{S}_w\coloneqq\bigcup\{F \subset Q \mid F \text{ is a facet of } Q \text{ with } w\in F \},
		\end{equation}
		we observe that $\rho_G (S_v) = S_v$ implies $\rho_G (\tilde{S}_w) = \tilde{S}_w$. Because $P$ is a local maximizer of the isotropic constant, Theorem \ref{thm_ccg} and Lemma \ref{lemma_pyramid_rs_decomposable} imply that $Q$ is RS-in\-de\-com\-po\-sable. Repeating the construction above for $Q$ and $w$ instead of $P$ and $v$ and %
		using that $Q$ is RS-indecomposable, 
		we arrive at the conclusion that $Q$ is a pyramid of the form $Q=\conv[Q_w \cup \{w\}]$ for some $(n-2)$-dimensional polytope $Q_w \subset H_v$. Since $Q$ is simplicial, it is a pyramid over a simplex and hence itself a simplex. Finally, since $P$ is a pyramid over $Q$, it is also a simplex.
	\end{proof}
	
	\begin{remark}
		We have already mentioned in Section \ref{sect_first_order_conditions} that \eqref{eq_spop_condition} is satisfied by a local maximizer $P$ of $K \mapsto L_K$ in the class of polytopes with at most $f_{n-1}(P)$ facets. Since the RS-movement $(P_t)_{t\in[-\varepsilon,\varepsilon]}$ from the proof of Theorem \ref{thm_simplicial_vertex} does not leave the class of polytopes with at most $f_{n-1}(P)$ facets, Theorem \ref{thm_simplicial_vertex} remains true if we interpret the expression ``locally maximizes'' as referring to local maximization within the class of polytopes with at most $f_{n-1}(P)$ facets.
	\end{remark}
	
	\section{The centrally symmetric case} \label{sect_the_centrally_symmetric_case}
	
	In this section, we prove Theorem \ref{thm_symmetric_simplicial_vertex}, which provides an analogue of Theorem \ref{thm_simplicial_vertex} for the centrally symmetric case. We follow Meckes \cite{mec05}, who adapted the technique of RS-movements to the centrally symmetric case.%
	
	\begin{definition}
		An RS-movement of a centrally symmetric convex body $K \subset \R^n$ with center of symmetry $c_K \in \R^n$ is called \textit{SRS-movement} if its speed function $\beta$ is \textit{odd}, i.e., if $\beta(x+c_K)=-\beta(-x+c_K)$ for all $x \in \pi_{v^\bot}(K-c_K)$. We say that a convex body $K \subset \R^n$ is \textit{SRS-indecomposable} if every SRS-movement $(K_t)_{t \in [-\varepsilon,\varepsilon]}$ with $K=K_0$ and $\varepsilon>0$ has an affine speed function.
	\end{definition}
	
	Let $(K_t)_{t \in [-\varepsilon,\varepsilon]}$ be an SRS-movement of a centrally symmetric convex body with $\varepsilon>0$. Observing that $K_t$ is centrally symmetric for every $t \in [-\varepsilon,\varepsilon]$, we obtain the following immediate corollary to Theorem \ref{thm_ccg}, which appears in \cite[Cor.~8]{mec05}.
	
	\begin{corollary} \label{cor_ccg_symmetric}
		Every local maximizer of the isotropic constant in the class of centrally symmetric convex bodies is SRS-indecomposable.
	\end{corollary}
	
	A key step towards Theorem \ref{thm_symmetric_simplicial_vertex} is the fact that \eqref{eq_spop_condition} holds for polytopal extremals in the class of centrally symmetric bodies, which is established by the following lemma.
	
	\begin{lemma}
		Let $P \subset \R^n$ be a centrally symmetric polytope that is a local maximizer of $K \mapsto L_K$ in the space of centrally symmetric convex bodies. Moreover, let $F \subset P$ be a facet and let $X$ be a random vector that is uniformly distributed on $F$. Then \eqref{eq_spop_condition} holds, i.e., we have $\E[\norm{X}^2 X] = (n+2) \E[X]$.
	\end{lemma}
	\begin{proof}
		Without loss of generality, we assume that $c_K=0$, i.e., that $P$ is $o$-symmetric. Let $f\colon \bd P \rightarrow R$ be an even function with $\supp f=  F \cup -F$ and whose restriction to $F$ is affine. Performing the construction from \cite[Sect.~3]{kip24}, we obtain a family of centrally symmetric polytopes $P_t$, $t \in [-1,1]$, such that the function $g \colon [-1,1] \rightarrow \R$, $t \mapsto L_{P_t}$ is differentiable at $t=0$. By assumption, $t=0$ is a local maximizer of $g$, and it follows that
		\begin{align}
			0=g'(0) &= \int_F (\norm{x}_2^2-n-2)\cdot f(x)\dint x+\int_{-F} (\norm{x}_2^2-n-2)\cdot f(x) \dint x\\&=2 \cdot \int_F (\norm{x}_2^2-n-2)\cdot f(x) \dint x.
		\end{align}
		In other words, we have \eqref{eq_spop_linear_condition} and hence also \eqref{eq_spop_condition}.
	\end{proof}
	
	Let $K \subset \R^n$ be a centrally symmetric $(n-1)$-dimensional convex body with center of symmetry $c \in K$. A \textit{bipyramid over $K$} is a convex body of the form
	\begin{equation}
		\tilde{K} = \conv[ (K-c) \cup \{p,-p\}]+c
	\end{equation}
	for some point $p \notin \aff (K-c)$. If $P \subset \R^n$ is a centrally symmetric polytope and $v \in P$ is a vertex, we say that $P$ is a bipyramid \textit{under $v$} if $P$ is a bipyramid over a hyperplane section $P \cap H_v$ with $v \notin H_v$. In the proof of Theorem \ref{thm_symmetric_simplicial_vertex}, we will use the fact that a centrally symmetric polytope that is a bipyramid under each of its vertices $v \in P$ is a cross-polytope. %
	To see this, let $P=-P$ be such a polytope and $v \in \vertices P$. If $P$ is a bipyramid over $Q \coloneqq P \cap H_v$, then $(\vertices P) \setminus \{v,-v\} = \vertices Q$ and $Q$ is also a bipyramid under each of its vertices. By induction over the dimension, $P$ is an iterated bipyramid over a point, i.e., a cross-polytope.

	\begin{proof}[Proof of Theorem \ref{thm_symmetric_simplicial_vertex}]
		Again, we assume without loss of generality that $P$ is isotropic. In particular, the center of symmetry of $P$ is the origin. Let $v \in P$ be a simplicial vertex. Because $P=-P$, the point $-v$ is also a simplicial vertex of $P$. We define $S_v$ and $S_{-v}$ as in Lemma \ref{lemma_local_symmetry} and observe that $S_v=-S_{-v}$. We also note that no facet contains both $v$ and $-v$, because otherwise the line segment $[-v,v]$ would be contained in that facet, in contradiction to the fact that $0\in [-v,v]$ lies in the interior of $P$. We define $\tilde{\beta}\colon \bd P \rightarrow \R$ by setting $\tilde{\beta}(v)=1$, $\tilde{\beta}(-v)=-1$ and $\tilde{\beta}(w)=0$ for all $w \in \vertices  P\setminus \{v,-v\}$, and interpolating affinely on each facet. As in the proof of Theorem \ref{thm_simplicial_vertex}, this definition is justified because every facet that contains $v$ or $-v$ is a simplex. 
		Let $G\subset P$ be a ridge that contains $v$. Then $-G$ contains $-v$ and we have $H \coloneqq \linh G = \linh (-G)$. By Lemma \ref{lemma_local_symmetry}, the sets $S_v$ and $S_{-v}$ are both symmetric with respect to $H$, which implies that $\tilde{\beta}$ factors through $\pi_H$ with, say, $\tilde{\beta} = \beta \circ \pi_H$. Let $u \in \sphere$ be a normal vector of $H$ and set
		\begin{equation}
			P_t \coloneqq \{x+t \beta(\pi_{H}(x))u \mid x \in P\} \quad \text{for } t \in \R.
		\end{equation}
		
		\begin{figure}[ht]
			\begin{subfigure}[b]{.33\linewidth}
				\centering
				\begin{tikzpicture}%
	[x={(0.562037cm, -0.330669cm)},
	y={(0.827112cm, 0.224669cm)},
	z={(0.000024cm, 0.916614cm)},
	scale=1.200000,
	back/.style={dotted, thin},
	edge/.style={color=black},
	facet/.style={fill=white,fill opacity=0.000000},
	vertex/.style={},
	back2/.style={dotted},
	edge2/.style={color=red,dashed},
	facet2/.style={fill=orange,fill opacity=0.20000},
	vertex2/.style={inner sep=1pt,circle,draw=darkgray!25!black,fill=darkgray!75!black,thick},
	back3/.style={dotted, thin},
	edge3/.style={color=blue,very thick},
	facet3/.style={fill=green,fill opacity=0.20000},
	vertex3/.style={inner sep=1pt,circle,draw=darkgray!25!black},
	vertex4/.style={inner sep=1pt,circle,draw=blue,fill=blue},
	facet4/.style={fill=blue,fill opacity=0.20000},]
	
	\coordinate (1.00000, -1.05000, 1.00000) at (1.00000, -1.05000, 1.00000);
	\coordinate (-1.00000, -1.05000, -1.00000) at (-1.00000, -1.05000, -1.00000);
	\coordinate (-1.00000, -0.95000, 1.00000) at (-1.00000, -0.95000, 1.00000);
	\coordinate (1.00000, -0.95000, -1.00000) at (1.00000, -0.95000, -1.00000);
	\coordinate (0.00000, -1.30000, 0.00000) at (0.00000, -1.30000, 0.00000);
	\coordinate (-1.00000, 1.05000, -1.00000) at (-1.00000, 1.05000, -1.00000);
	\coordinate (1.00000, 1.05000, 1.00000) at (1.00000, 1.05000, 1.00000);
	\coordinate (1.00000, 0.95000, -1.00000) at (1.00000, 0.95000, -1.00000);
	\coordinate (-1.00000, 0.95000, 1.00000) at (-1.00000, 0.95000, 1.00000);
	\coordinate (0.00000, 1.30000, 0.00000) at (0.00000, 1.30000, 0.00000);
	\draw[edge,back] (-1.00000, -1.05000, -1.00000) -- (-1.00000, 1.05000, -1.00000);
	\draw[edge,back] (-1.00000, 1.05000, -1.00000) -- (1.00000, 0.95000, -1.00000);
	\draw[edge,back] (-1.00000, 1.05000, -1.00000) -- (-1.00000, 0.95000, 1.00000);
	\draw[edge,back] (-1.00000, 1.05000, -1.00000) -- (0.00000, 1.30000, 0.00000);
	\draw[edge,back] (1.00000, 1.05000, 1.00000) -- (0.00000, 1.30000, 0.00000);
	\draw[edge,back] (1.00000, 0.95000, -1.00000) -- (0.00000, 1.30000, 0.00000);
	\draw[edge,back] (-1.00000, 0.95000, 1.00000) -- (0.00000, 1.30000, 0.00000);
	\node[vertex] at (-1.00000, 1.05000, -1.00000)     {};
	\node[vertex] at (0.00000, 1.30000, 0.00000)     {};
	\fill[facet4] (1.00000, 1.30000, -1.00000) -- (-1.00000, 1.30000, 1.00000) -- (-1.00000, -1.3000, 1.0000) -- (1.00000, -1.3000, -1.00000)  -- cycle {};
	\draw[edge] (1.00000, -1.05000, 1.00000) -- (-1.00000, -0.95000, 1.00000);
	\draw[edge] (1.00000, -1.05000, 1.00000) -- (1.00000, -0.95000, -1.00000);
	\draw[edge] (1.00000, -1.05000, 1.00000) -- (0.00000, -1.30000, 0.00000);
	\draw[edge] (1.00000, -1.05000, 1.00000) -- (1.00000, 1.05000, 1.00000);
	\draw[edge] (-1.00000, -1.05000, -1.00000) -- (-1.00000, -0.95000, 1.00000);
	\draw[edge] (-1.00000, -1.05000, -1.00000) -- (1.00000, -0.95000, -1.00000);
	\draw[edge] (-1.00000, -1.05000, -1.00000) -- (0.00000, -1.30000, 0.00000);
	\draw[edge] (-1.00000, -0.95000, 1.00000) -- (0.00000, -1.30000, 0.00000);
	\draw[edge] (-1.00000, -0.95000, 1.00000) -- (-1.00000, 0.95000, 1.00000);
	\draw[edge3] (1.00000, -0.95000, -1.00000) -- (0.00000, -1.30000, 0.00000);
	\draw[edge] (1.00000, -0.95000, -1.00000) -- (1.00000, 0.95000, -1.00000);
	\draw[edge] (1.00000, 1.05000, 1.00000) -- (1.00000, 0.95000, -1.00000);
	\draw[edge] (1.00000, 1.05000, 1.00000) -- (-1.00000, 0.95000, 1.00000);

	\draw[edge2] (-1.00000, 1.30000, -1.00000) -- (1.00000, 1.30000, 1.00000);
	\draw[edge2] (1.00000, -1.30000, 1.00000) -- (-1.00000, -1.30000, -1.00000);
	\node[vertex] at (1.00000, -1.05000, 1.00000)     {};
	\node[vertex] at (-1.00000, -1.05000, -1.00000)     {};
	\node[vertex] at (-1.00000, -0.95000, 1.00000)     {};
	\node[vertex] at (1.00000, -0.95000, -1.00000)     {};
	\node[vertex2] at (0.00000, -1.30000, 0.00000)     {};
	\node[vertex3] at (1.00000, -1.30000, 1.00000)     {};
	\node[vertex3] at (-1.00000, -1.30000, -1.00000)     {};
	\node[vertex4] at (0.00000, 0.0000, 0.00000)     {};
	\node[vertex] at (1.00000, 1.05000, 1.00000)     {};
	\node[vertex] at (1.00000, 0.95000, -1.00000)     {};
	\node[vertex] at (-1.00000, 0.95000, 1.00000)     {};
	
	\node[vertex2] at (0.00000, 1.30000, 0.00000)     {};
	\node[vertex3] at (-1.00000, 1.30000, -1.00000)     {};
	\node[vertex3] at (1.00000, 1.30000, 1.00000)     {};
\end{tikzpicture}
				
				\caption{initial polytope $P_0$}
			\end{subfigure}
			\begin{subfigure}[b]{.32\linewidth}
				\centering
				\begin{tikzpicture}%
	[x={(0.562037cm, -0.330669cm)},
	y={(0.827112cm, 0.224669cm)},
	z={(0.000024cm, 0.916614cm)},
	scale=1.200000,
	back/.style={dotted, thin},
	edge/.style={color=black},
	facet/.style={fill=white,fill opacity=0.000000},
	vertex/.style={},
	back2/.style={dotted},
	edge2/.style={color=red,dashed},
	facet2/.style={fill=Emerald,fill opacity=0.20000},
	vertex2/.style={inner sep=1pt,circle,draw=darkgray!25!black,fill=darkgray!75!black,thick},
	back3/.style={dotted, thin},
	edge3/.style={color=black,very thick},
	facet3/.style={fill=Emerald,fill opacity=0.10000},
	vertex3/.style={inner sep=1pt,circle,draw=darkgray!25!black},
	facet4/.style={fill=blue,fill opacity=0.20000},]
	
	\coordinate (1.00000, -1.05000, 1.00000) at (1.00000, -1.05000, 1.00000);
	\coordinate (-1.00000, -1.05000, -1.00000) at (-1.00000, -1.05000, -1.00000);
	\coordinate (-1.00000, -0.95000, 1.00000) at (-1.00000, -0.95000, 1.00000);
	\coordinate (1.00000, -0.95000, -1.00000) at (1.00000, -0.95000, -1.00000);
	\coordinate (-1.00000, 1.05000, -1.00000) at (-1.00000, 1.05000, -1.00000);
	\coordinate (1.00000, 1.05000, 1.00000) at (1.00000, 1.05000, 1.00000);
	\coordinate (1.00000, 0.95000, -1.00000) at (1.00000, 0.95000, -1.00000);
	\coordinate (-1.00000, 0.95000, 1.00000) at (-1.00000, 0.95000, 1.00000);
	\coordinate (-0.40000, -1.30000, -0.40000) at (-0.40000, -1.30000, -0.40000);
	\coordinate (0.40000, 1.30000, 0.40000) at (0.40000, 1.30000, 0.40000);
	\draw[edge,back] (-1.00000, -1.05000, -1.00000) -- (-1.00000, 1.05000, -1.00000);
	\draw[edge,back] (-1.00000, 1.05000, -1.00000) -- (1.00000, 0.95000, -1.00000);
	\draw[edge,back] (-1.00000, 1.05000, -1.00000) -- (-1.00000, 0.95000, 1.00000);
	\draw[edge,back] (-1.00000, 1.05000, -1.00000) -- (0.40000, 1.30000, 0.40000);
	\draw[edge,back] (1.00000, 1.05000, 1.00000) -- (0.40000, 1.30000, 0.40000);
	\draw[edge,back] (1.00000, 0.95000, -1.00000) -- (0.40000, 1.30000, 0.40000);
	\draw[edge,back] (-1.00000, 0.95000, 1.00000) -- (0.40000, 1.30000, 0.40000);
	\node[vertex] at (-1.00000, 1.05000, -1.00000)     {};
	\node[vertex] at (0.40000, 1.30000, 0.40000)     {};
	\fill[facet2] (-0.40000, -1.30000, -0.40000) -- (1.00000, -1.05000, 1.00000) -- (-1.00000, -0.95000, 1.00000) -- cycle {};
	\fill[facet2] (-0.40000, -1.30000, -0.40000) -- (-1.00000, -1.05000, -1.00000) -- (-1.00000, -0.95000, 1.00000) -- cycle {};
	\fill[facet2] (-0.40000, -1.30000, -0.40000) -- (-1.00000, -1.05000, -1.00000) -- (1.00000, -0.95000, -1.00000) -- cycle {};
	\fill[facet2] (-0.40000, -1.30000, -0.40000) -- (1.00000, -1.05000, 1.00000) -- (1.00000, -0.95000, -1.00000) -- cycle {};
	
	\fill[facet3] (0.40000, 1.30000, 0.40000) -- (-1.00000, 1.05000, -1.00000) -- (1.00000, 0.95000, -1.00000) -- cycle {};
	\fill[facet3] (0.40000, 1.30000, 0.40000) -- (1.00000, 1.05000, 1.00000) -- (1.00000, 0.95000, -1.00000) -- cycle {};
	\fill[facet3] (0.40000, 1.30000, 0.40000) -- (1.00000, 1.05000, 1.00000) -- (-1.00000, 0.95000, 1.00000) -- cycle {};
	\fill[facet3] (0.40000, 1.30000, 0.40000) -- (-1.00000, 1.05000, -1.00000) -- (-1.00000, 0.95000, 1.00000) -- cycle {};
	
	\draw[edge] (1.00000, -1.05000, 1.00000) -- (-1.00000, -0.95000, 1.00000);
	\draw[edge] (1.00000, -1.05000, 1.00000) -- (1.00000, -0.95000, -1.00000);
	\draw[edge] (1.00000, -1.05000, 1.00000) -- (1.00000, 1.05000, 1.00000);
	\draw[edge] (1.00000, -1.05000, 1.00000) -- (-0.40000, -1.30000, -0.40000);
	\draw[edge] (-1.00000, -1.05000, -1.00000) -- (-1.00000, -0.95000, 1.00000);
	\draw[edge] (-1.00000, -1.05000, -1.00000) -- (1.00000, -0.95000, -1.00000);
	\draw[edge] (-1.00000, -1.05000, -1.00000) -- (-0.40000, -1.30000, -0.40000);
	\draw[edge] (-1.00000, -0.95000, 1.00000) -- (-1.00000, 0.95000, 1.00000);
	\draw[edge] (-1.00000, -0.95000, 1.00000) -- (-0.40000, -1.30000, -0.40000);
	\draw[edge] (1.00000, -0.95000, -1.00000) -- (1.00000, 0.95000, -1.00000);
	\draw[edge] (1.00000, -0.95000, -1.00000) -- (-0.40000, -1.30000, -0.40000);
	\draw[edge] (1.00000, 1.05000, 1.00000) -- (1.00000, 0.95000, -1.00000);
	\draw[edge] (1.00000, 1.05000, 1.00000) -- (-1.00000, 0.95000, 1.00000);
	
	\draw[edge2] (-1.00000, 1.30000, -1.00000) -- (1.00000, 1.30000, 1.00000);
	\draw[edge2] (1.00000, -1.30000, 1.00000) -- (-1.00000, -1.30000, -1.00000);
	\node[vertex] at (1.00000, -1.05000, 1.00000)     {};
	\node[vertex] at (-1.00000, -1.05000, -1.00000)     {};
	\node[vertex] at (-1.00000, -0.95000, 1.00000)     {};
	\node[vertex] at (1.00000, -0.95000, -1.00000)     {};
	\node[vertex] at (1.00000, 1.05000, 1.00000)     {};
	\node[vertex] at (1.00000, 0.95000, -1.00000)     {};
	\node[vertex] at (-1.00000, 0.95000, 1.00000)     {};
	\node[vertex3] at (1.00000, -1.30000, 1.00000)     {};
	\node[vertex3] at (-1.00000, -1.30000, -1.00000)     {};
	\node[vertex2] at (-0.40000, -1.30000, -0.40000)     {};
	
	\node[vertex3] at (-1.00000, 1.30000, -1.00000)     {};
	\node[vertex3] at (1.00000, 1.30000, 1.00000)     {};
	\node[vertex2] at (0.40000, 1.30000, 0.40000)     {};
\end{tikzpicture}
				
				\caption{$\pm v$ moves in direction $\pm u$}
			\end{subfigure}
			\begin{subfigure}[b]{.32\linewidth}
				\centering
				\begin{tikzpicture}%
	[x={(0.562037cm, -0.330669cm)},
	y={(0.827112cm, 0.224669cm)},
	z={(0.000024cm, 0.916614cm)},
	scale=1.200000,
	back/.style={dotted, thin},
	edge/.style={color=black},
	facet/.style={fill=white,fill opacity=0.000000},
	vertex/.style={},
	back2/.style={dotted},
	edge2/.style={color=red,dashed},
	facet2/.style={fill=Emerald,fill opacity=0.20000},
	vertex2/.style={inner sep=1pt,circle,draw=darkgray!25!black,fill=darkgray!75!black,thick},
	back3/.style={dotted, thin},
	edge3/.style={color=black,very thick},
	facet3/.style={fill=Emerald,fill opacity=0.10000},
	vertex3/.style={inner sep=1pt,circle,draw=darkgray!25!black},
	facet4/.style={fill=blue,fill opacity=0.20000},]
	
	\coordinate (1.00000, -1.05000, 1.00000) at (1.00000, -1.05000, 1.00000);
	\coordinate (-1.00000, -1.30000, -1.00000) at (-1.00000, -1.30000, -1.00000);
	\coordinate (-1.00000, -0.95000, 1.00000) at (-1.00000, -0.95000, 1.00000);
	\coordinate (1.00000, -0.95000, -1.00000) at (1.00000, -0.95000, -1.00000);
	\coordinate (1.00000, 0.95000, -1.00000) at (1.00000, 0.95000, -1.00000);
	\coordinate (-1.00000, 0.95000, 1.00000) at (-1.00000, 0.95000, 1.00000);
	\coordinate (-1.00000, 1.05000, -1.00000) at (-1.00000, 1.05000, -1.00000);
	\coordinate (1.00000, 1.30000, 1.00000) at (1.00000, 1.30000, 1.00000);
	\draw[edge,back] (-1.00000, -1.30000, -1.00000) -- (-1.00000, 1.05000, -1.00000);
	\draw[edge,back] (1.00000, 0.95000, -1.00000) -- (-1.00000, 1.05000, -1.00000);
	\draw[edge,back] (-1.00000, 0.95000, 1.00000) -- (-1.00000, 1.05000, -1.00000);
	\draw[edge,back] (-1.00000, 1.05000, -1.00000) -- (1.00000, 1.30000, 1.00000);
	\node[vertex] at (-1.00000, 1.05000, -1.00000)     {};
	\fill[facet2] (-1.00000, -0.95000, 1.00000) -- (1.00000, -1.05000, 1.00000) -- (-1.00000, -1.30000, -1.00000) -- cycle {};
	\fill[facet2] (1.00000, -0.95000, -1.00000) -- (1.00000, -1.05000, 1.00000) -- (-1.00000, -1.30000, -1.00000) -- cycle {};
	
	\fill[facet3] (1.00000, 0.95000, -1.00000) -- (-1.00000, 1.05000, -1.00000) -- (1.00000, 1.30000, 1.00000) -- cycle {};
	\fill[facet3] (-1.00000, 0.95000, 1.00000) -- (-1.00000, 1.05000, -1.00000) -- (1.00000, 1.30000, 1.00000) -- cycle {};
	\draw[edge] (1.00000, -1.05000, 1.00000) -- (-1.00000, -1.30000, -1.00000);
	\draw[edge] (1.00000, -1.05000, 1.00000) -- (-1.00000, -0.95000, 1.00000);
	\draw[edge] (1.00000, -1.05000, 1.00000) -- (1.00000, -0.95000, -1.00000);
	\draw[edge] (1.00000, -1.05000, 1.00000) -- (1.00000, 1.30000, 1.00000);
	\draw[edge] (-1.00000, -1.30000, -1.00000) -- (-1.00000, -0.95000, 1.00000);
	\draw[edge] (-1.00000, -1.30000, -1.00000) -- (1.00000, -0.95000, -1.00000);
	\draw[edge] (-1.00000, -0.95000, 1.00000) -- (-1.00000, 0.95000, 1.00000);
	\draw[edge] (1.00000, -0.95000, -1.00000) -- (1.00000, 0.95000, -1.00000);
	\draw[edge] (1.00000, 0.95000, -1.00000) -- (1.00000, 1.30000, 1.00000);
	\draw[edge] (-1.00000, 0.95000, 1.00000) -- (1.00000, 1.30000, 1.00000);
	
	\draw[edge2] (1.00000, -1.30000, 1.00000) -- (-1.00000, -1.30000, -1.00000);
	\draw[edge2] (-1.00000, 1.30000, -1.00000) -- (1.00000, 1.30000, 1.00000);
	\node[vertex] at (1.00000, -1.05000, 1.00000)     {};
	\node[vertex] at (-1.00000, -1.30000, -1.00000)     {};
	\node[vertex] at (-1.00000, -0.95000, 1.00000)     {};
	\node[vertex] at (1.00000, -0.95000, -1.00000)     {};
	\node[vertex] at (1.00000, 0.95000, -1.00000)     {};
	\node[vertex] at (-1.00000, 0.95000, 1.00000)     {};
	\node[vertex] at (1.00000, 1.30000, 1.00000)     {};
	
	\node[vertex3] at (1.00000, -1.30000, 1.00000)     {};
	\node[vertex2] at (-1.00000, -1.30000, -1.00000)     {};
	\node[vertex3] at (-1.00000, 1.30000, -1.00000)     {};
	\node[vertex2] at (1.00000, 1.30000, 1.00000)     {};
\end{tikzpicture}
				
				\caption{end of movement $P_\varepsilon$}
			\end{subfigure}
			\caption{The SRS-movement $(P_t)_{i \in [-\varepsilon,\varepsilon]}$. The vertices $v$ and $-v$ move orthogonally to the linear hyperplane spanned by the ridge $G$ and the origin (marked in blue).}
		\end{figure}
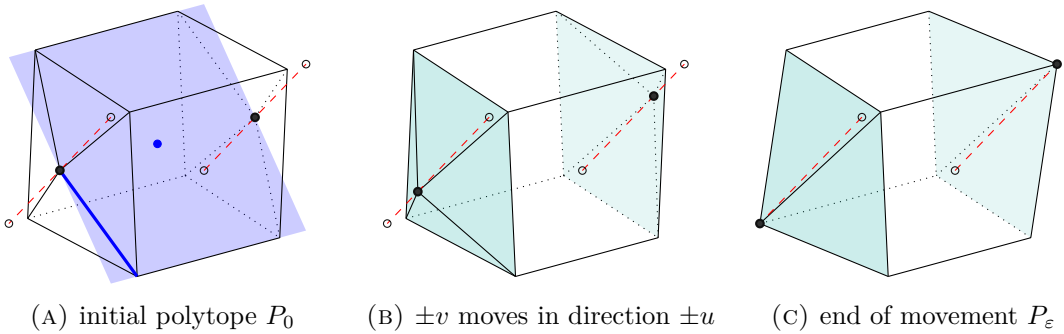
	
		As in the proof of Theorem \ref{thm_simplicial_vertex}, we observe that there exists a $\varepsilon>0$ such that $(P_t)_{[-\varepsilon,\varepsilon]}$ is an SRS-movement of $P$. By Corollary \ref{cor_ccg_symmetric}, it follows that $\beta$ is linear. This implies that $\pi_H(S_v \cup S_{-v}) = \pi_H(P)$ and that there exists a hyperplane $H_v$ with $\vertices  P\setminus \{v,-v\} \subset H_v$. In other words, $P$ is a bipyramid under $v$. But this implies that $S_v \cup S_{-v} = \bd P$ and hence that $P$ is simplicial. Repeating the argument above for an arbitrary vertex $w \in P$, we obtain that $P$ is a bipyramid under $w$ for every $w \in \vertices P$ and hence a cross-polytope.
	\end{proof}
	
	\section{Cubical zones of zonotopes} \label{sect_cubical_zones_of_zonotopes}
	
	In the preceding two sections, we have considered affine reflectors that are contained in two simplices. Another interesting example is given by the case where the facets $F_1$ and $F_2$ are prisms over $F_1 \cap F_2$, which implies again that $F_1 \cap F_2$ is an affine reflector. The analogue of a simplicial polytope in this setting is a polytope whose facets are all parallelepipeds. As mentioned in the introduction, such a polytope has to be a zonotope. For basic facts about zonotopes, we refer to \cite[Ch.~7]{zie07} and the references therein. 
	
	We note that Theorem \ref{thm_linear_reflector_simplex} implies that a zonotope $Z$ whose facets are all parallelepipeds cannot be a local maximizer of the isotropic constant. However, since $Z$ is necessarily centrally symmetric, it is natural to ask what can be said about $Z$ if $Z$ is assumed to be a local maximizer of $K \mapsto L_K$ in the class of centrally symmetric convex bodies. In this section, we prove Theorem \ref{thm_cubical_zone}, which implies that such a zonotope $Z$ has to be a cube. As in the preceding sections, it turns out to be sufficient that the boundary of $Z$ contains enough affine reflectors to enforce a local hyperplane symmetry, which leads us to consider zonotopes with at least one cubical zone.
	
	Let $Z \subset \R^n$ be a full-dimensional zonotope. By definition, $Z$ can be written as
	\begin{equation}
		Z = c+\sum_{i=1}^m [-z_i,z_i]
	\end{equation}
	for certain vectors $z_1,\dots,z_m \in \R^n \setminus \{0\}$, $m\geq n$, and a point $c \in \R^n$. We identify the vectors $z_1,\dots,z_m$ with the line segments $[-z_1,z_1],\dots,[-z_m,z_m]$ and call them \textit{generators} of $Z$. A set of generators $z_1,\dots,z_m$ is said to be \textit{irredundant} if the vectors are pairwise non-parallel. In this case, there is a one-to-one correspondence between the generators $z_i$ and the \textit{zones} of $Z$, which are given by
	\begin{equation}
		\mathcal{Z}_i \coloneqq \{F \subset Z \mid F \text{ proper face of } Z \text{ with $[-z_i,z_i]$ as a Minkowski summand}\}.
	\end{equation}
	We fix an index $i \in [m]$ and consider the orthogonal projection $\pi_{z_i^\bot}(Z)$ of $Z$ onto the hyperplane $z_i^\bot$. Then
	$\pi_{z_i^\bot}(Z)$ is again a zonotope, whose proper faces are in one-to-one correspondence with the elements of the zone $\mathcal{Z}_i$. If we assume that $\mathcal{Z}_i$ is cubical, i.e., that all elements of $\mathcal{Z}_i$ are parallelepipeds, then every ridge in $\mathcal{Z}_i$ is an affine reflector. Adapting the proof of Lemma \ref{lemma_symmetric_wrt_hyperplane} once again, we obtain the following lemma.	
	
	\begin{lemma} \label{lemma_zone_hyperplane_symmetry}
		Let $Z$ be an isotropic zonotope that satisfies \eqref{eq_spop_condition}. If $Z$ has a cubical zone $\mathcal{Z}$, then $\bigcup \mathcal{Z}$ is symmetric with respect to the linear hyperplane spanned by any ridge $G \in \mathcal{Z}$.
	\end{lemma}
	
	\begin{proof}
		Let $z_1,\dots,z_m$ be irredundant generators of $Z$ and assume without loss of generality that $\mathcal{Z}=\mathcal{Z}_1$. Because $Z$ satisfies \eqref{eq_spop_condition}, Lemma \ref{lemma_euclidean_reflector} implies that every ridge $G \in \mathcal{Z}$ is a Euclidean reflector. We fix a ridge $G = F_1 \cap F_{-1} \in \mathcal{Z}$ and consider an arbitrary facet $F \in \mathcal{Z}$. Since the facets in $\mathcal{Z}$ are in one-to-one correspondence with the facets of $\pi_{z_1^\bot}(Z)$, there exists a sequence of facets $F_1,F_2,\dots,F_\ell=F$ in $\mathcal{Z}$ such that $F_i$ and $F_{i+1}$ are adjacent for all $i \in [\ell-1]$. Applying the same argument as in the proof of Lemma \ref{lemma_zone_hyperplane_symmetry}, we conclude that $\rho_G(F) \in \mathcal{Z}$ is a facet of $Z$.
	\end{proof}
	
	Let $K \subset \R^n$ be a centrally symmetric convex body. If there exist a centrally symmetric $(n-1)$-dimensional convex body $\tilde{K}$ and a point $v \in \R^n \setminus (\linh \tilde{K})$ such that $K$ is given by
	\begin{equation}
		K = \tilde{K}+[-v,v],
	\end{equation}
	we say that $K$ is a \textit{prism over $\tilde{K}$} and a \textit{prism under $[-v,v]$}. A simple computation shows that the isotropic constant of an isotropic prism $K$ over $\tilde{K}$ is given by
	\begin{equation} \label{eq_L_K_prism_over_Q}
		L_K^n = \frac{1}{\vol K} = \frac{1}{\sqrt{12} \vol \tilde{K}} = \frac{1}{\sqrt{12}} L_{\tilde{K}}^{n-1}.
	\end{equation}
	The key idea of the proof of Theorem \ref{thm_cubical_zone} is that we can use \eqref{eq_spop_condition} and Corollary \ref{cor_ccg_symmetric} to show that $Z$ must be an iterated prism over a line segment.	
	\begin{proof}[Proof of Theorem \ref{thm_cubical_zone}]
		Without loss of generality, we assume that $Z$ is isotropic and that $Z$ is given by irredundant generators $z_1,\dots,z_m \in \R^n$.
		Let $\mathcal{Z}$ be a cubical zone of $Z$, say, $\mathcal{Z}=\mathcal{Z}_1$. Let $L \coloneqq \linh \{z_1\}$. Every vertex $v$ of $Z$ is an extreme point of the line segment $I_v \coloneqq (L+v) \cap Z$. We set
		\begin{equation}
			\tilde{\beta}(v)=\begin{cases}
				1 & \text{if } \scpr{z_1,v} = \max_{x\in I_v} \scpr{z_1,x},\\
				-1 & \text{if } \scpr{z_1,v} = \min_{x\in I_v} \scpr{z_1,x}
			\end{cases}
		\end{equation}
		and extend $\tilde{\beta}$ to a function $\bd Z \rightarrow \R$ by interpolating affinely on each facet. This is well-defined because every facet $F \subset Z$ that contains two vertices $v,w$ with $\tilde{\beta}(v)=1$, $\tilde{\beta}(w)=-1$ has the generator $[-z_1,z_1]$ as a Minkowski summand and is hence, by assumption, a parallelepiped and in particular a prism under $[-z_1,z_1]$. On each facet outside the zone $\mathcal{Z}$, $\tilde{\beta}$ is constant with value $1$ or $-1$.
		
		Let $G \in \mathcal{Z}$ be an arbitrary ridge. By Lemma \ref{lemma_zone_hyperplane_symmetry}, $\bigcup\mathcal{Z}$ is symmetric with respect to $H\coloneqq \linh G$. Therefore, the function $\tilde{\beta}$ factors through $\pi_H$, i.e., there exists $\beta \colon \pi_H(Z) \rightarrow \R$ with $\tilde{\beta} = \beta \circ \pi_H$. Moreover, if $u\in \sphere$ is a normal vector of $H$, there exists an $\varepsilon>0$ such that
		\begin{equation}
			Z_t \coloneqq \{x+t \beta(\pi_{H}(x))u \mid x \in P\}=[-(z_1+tu),z_1+tu] +\sum_{i=2}^m [-z_i,z_i] \quad \text{for } t \in [-\varepsilon,\varepsilon].
		\end{equation}
		By Corollary \ref{cor_ccg_symmetric}, $\beta$ must be linear. It follows that $Z$ is a prism under $[-z_1,z_1]$.
		
		Since $Z$ is a local maximizer of the isotropic constant in the class of $n$-dimensional centrally symmetric convex bodies, \eqref{eq_L_K_prism_over_Q} implies that $\pi_{z_1^\bot}(Z)$ is a local maximizer in the class of $(n-1)$-dimensional centrally symmetric convex bodies. Moreover, since the facets of $\pi_{z_1^\bot}(Z)$ are given by $\{\pi_{z_1^\bot}(F) \mid F \in \mathcal{Z} \text{ is a facet}\}$
		and $\mathcal{Z}$ is assumed to be cubical, it follows that all facets of $\pi_{z_1^\bot}(Z)$ are parallelepipeds. In particular, $\pi_{z_1^\bot}(Z)$ has a cubical zone. By induction over the dimension, $Z$ is an iterated prism over a line segment and hence a cube.
	\end{proof}

	\section{Zonotopes with at most $n+1$ generators} \label{sect_zonotopes_n+1_generators}
	
	In the preceding sections, we obtained structural information about local maximizers of the isotropic constant under various assumptions concerning the presence of affine reflectors. However, the results above do not yield upper bounds on the isotropic constant (even for the respective subclasses of the space of convex bodies) because it is a priori unclear whether local maximizers with the respective properties exist. This is due to two independent reasons:
	\begin{enumerate}
		\item the classes under consideration (polytopes whose ridges are all affine reflectors, polytopes with a simplicial vertex, zonotopes with a cubical zone) are not closed with respect to the Banach-Mazur distance;
		\item the first-order condition \eqref{eq_spop_condition} is only valid for a local maximizer $P$ in the class of polytopes with at most $f_{n-1}(P)$ facets.
	\end{enumerate}
	In this section, we study a closed subclass of the Banach-Mazur compactum, namely, zonotopes with at most $n+1$ generators. We will establish in Lemma \ref{lemma_n+1_zonotopes_compactness} below that this class can be continuously parameterized by a compact set in $\R^n$, implying that the isotropic constant attains its minimum and its maximum in this class. The maximizers and the minimizers can then be identified by showing that almost all zonotopes with at most $n+1$ generators are SRS-decomposable, with SRS-movements that stay in the class of zonotopes with at most $n+1$ generators. In fact, we will see that any zonotope with exactly $n+1$ generators is symmetric enough to admit a non-trivial SRS-movement via a reduction to the planar case. In particular, no first-order conditions are needed.

	Throughout this section, all zonotopes are assumed to be full-dimensional. We denote the \textit{standard unit cube} in $\R^n$ by
	\begin{equation}
		C_n \coloneqq \sum_{i=1}^n [-e_i,e_i] \subset \R^n.
	\end{equation}
	
	\begin{proposition} \label{prop_n+1_zonotopes_upper_bound}
		Let $Z \subset \R^n$ be a zonotope with at most $n+1$ generators. Then
		\begin{equation}
			L_Z \leq L_{C_n}
		\end{equation}
		with equality if and only if $Z$ is affinely equivalent to $C_n$.
	\end{proposition}
	
	Since the planar case of \eqref{eq_symmetric_strong_isotropic_constant_conjecture} is known to be true \cite[Cor.~10]{mec05}, we know in particular that Proposition \ref{prop_n+1_zonotopes_upper_bound} holds in $\R^2$. Nevertheless, we start with a lemma that addresses the planar case of Proposition \ref{prop_n+1_zonotopes_upper_bound}. Here and in the following, $y \geq x$ for two vectors $x,y \in \R^n$ means that the inequality holds component-wise, i.e., we have $y_i \geq x_i$ for all $i \in [n]$. Similarly, $y > x$ means that we have $y_i > x_i$ for all $i \in [n]$.
	
	\begin{lemma} \label{lemma_2d}
		Let $x, y \in \R^2$ with $y \geq x$. Let $s_x, s_y \in \R$. Then
		\begin{equation}
			Z_t \coloneqq \left[x+ts_x\begin{bmatrix}
				1\\-1
			\end{bmatrix},y+ts_y\begin{bmatrix}
			1\\-1
		\end{bmatrix}\right]+C_2
		\end{equation}
		defines an RS-movement of the zonotope $Z_0=[x,y]+C_2$ on the interval
		\begin{equation}
			I \coloneqq \left\{ t \in \R \MID x+ts_x\begin{bmatrix}
				1\\-1
			\end{bmatrix} \leq y+ ts_y\begin{bmatrix}
			1\\-1
		\end{bmatrix}\right\}.
		\end{equation}
		If $x \neq y$ and $s_x \neq s_y$, then the corresponding speed function is not affine.
	\end{lemma}
	
	\begin{proof}
		By definition, $Z_t$ is a zonotope for all $t \in I$. We have to show that there exists a speed function $\beta$ %
		such that $Z_t$ can be written as in \eqref{eq_def_rs_movement}. Let $f \colon \R^2 \rightarrow \R$ be the linear functional given by $z \mapsto z_1+z_2$. We define three polyhedra
		\begin{align}
			P_- \coloneqq f^{-1}((-\infty,x_1+x_2]), \quad P_0 \coloneqq f^{-1}([x_1+x_2,y_1+y_2]), \quad P_+ \coloneqq f^{-1}([y_1+y_2,\infty))
		\end{align}
		and observe that
		\begin{equation}
			Z_t \coloneqq \conv[(Z_t \cap P_-) \cup (Z_t \cap P_+)] \quad \text{for } t\in I.
		\end{equation}
		
		\begin{figure}[ht]
			\begin{subfigure}[b]{.33\linewidth}
				\centering
				\begin{tikzpicture}[scale=0.6,
	back/.style={dotted, thin},
	edge/.style={color=black},
	facet/.style={fill=Emerald,fill opacity=0.200000},
	vertex/.style={},
	back2/.style={dotted},
	edge2/.style={color=red,dashed},
	facet2/.style={fill=Emerald,fill opacity=0.20000},
	vertex2/.style={inner sep=1pt,circle,draw=darkgray!25!black,fill=darkgray!75!black,thick},
	back3/.style={dotted, thin},
	edge3/.style={color=black,very thick},
	facet3/.style={fill=Emerald,fill opacity=0.10000},
	vertex3/.style={inner sep=1pt,circle,draw=darkgray!25!black},
	facet4/.style={fill=blue,fill opacity=0.20000},]
	\def\x{1.5};
	\def\y{1};
	
	\draw[back] (-1, -\x-\y-1) to (1,-\x-\y-1) to (1,\x+\y+1) to (-1,\x+\y+1) to (-1, -\x-\y-1);
	
	\draw[back] (-\x-\y-1, 1) to (-\x-\y-1,-1) to (\x+\y+1,-1) to (\x+\y+1,1) to (-\x-\y-1, 1);

	\fill[facet] (-\x-1,-\y+1) to (-\x-1,-\y-1) to (-\x+1,-\y-1) to (\x+1,\y-1) to (\x+1,\y+1) to (\x-1, \y+1) to (-\x-1,-\y+1);
	\draw[edge] (-\x-1,-\y+1) to (-\x-1,-\y-1) to (-\x+1,-\y-1) to (\x+1,\y-1) to (\x+1,\y+1) to (\x-1, \y+1) to (-\x-1,-\y+1);
	 
	\draw[edge2] (-1, -\x-\y-1) to (-\x-\y-1,-1);
	\draw[edge2] (1, -\x-\y-1) to (-\x-\y-1,+1);
	\draw[edge2] (1, \x+\y+1) to (\x+\y+1,1);
	\draw[edge2] (-1, \x+\y+1) to (\x+\y+1,-1);
	
	\node[vertex2] at (\x,\y) {};
	\node at (\x+0.4,\y+0.4) {$y$};
	\node[vertex2] at (-\x,-\y) {};
	\node at (-\x-0.4,-\y-0.4) {$x$};
	
\end{tikzpicture}
				
				\caption{initial zonotope $Z_0$}
			\end{subfigure}
			\begin{subfigure}[b]{.32\linewidth}
				\centering
				\begin{tikzpicture}%
	[scale=0.6,
	back/.style={dotted, thin},
	edge/.style={color=black},
	facet/.style={fill=Emerald,fill opacity=0.200000},
	vertex/.style={},
	back2/.style={dotted},
	edge2/.style={color=red,dashed},
	facet2/.style={fill=Emerald,fill opacity=0.20000},
	vertex2/.style={inner sep=1pt,circle,draw=darkgray!25!black,fill=darkgray!75!black,thick},
	back3/.style={dotted, thin},
	edge3/.style={color=black,very thick},
	facet3/.style={fill=Emerald,fill opacity=0.10000},
	vertex3/.style={inner sep=1pt,circle,draw=darkgray!25!black},
	facet4/.style={fill=blue,fill opacity=0.20000},]
	\def\x{1.9};
	\def\y{0.6};
	
	\draw[back] (-1, -\x-\y-1) to (1,-\x-\y-1) to (1,\x+\y+1) to (-1,\x+\y+1) to (-1, -\x-\y-1);
	
	\draw[back] (-\x-\y-1, 1) to (-\x-\y-1,-1) to (\x+\y+1,-1) to (\x+\y+1,1) to (-\x-\y-1, 1);

	\fill[facet] (-\x-1,-\y+1) to (-\x-1,-\y-1) to (-\x+1,-\y-1) to (\x+1,\y-1) to (\x+1,\y+1) to (\x-1, \y+1) to (-\x-1,-\y+1);
	\draw[edge] (-\x-1,-\y+1) to (-\x-1,-\y-1) to (-\x+1,-\y-1) to (\x+1,\y-1) to (\x+1,\y+1) to (\x-1, \y+1) to (-\x-1,-\y+1);
	
	\draw[edge2] (-1, -\x-\y-1) to (-\x-\y-1,-1);
	\draw[edge2] (1, -\x-\y-1) to (-\x-\y-1,+1);
	\draw[edge2] (1, \x+\y+1) to (\x+\y+1,1);
	\draw[edge2] (-1, \x+\y+1) to (\x+\y+1,-1);
	
	\node[vertex2] at (\x,\y) {};
	\node at (\x+0.4,\y+0.4) {$y$};
	\node[vertex2] at (-\x,-\y) {};
	\node at (-\x-0.4,-\y-0.4) {$x$};
\end{tikzpicture}
				
				\caption{movement along $\pm (e_1-e_2)$}
			\end{subfigure}
			\begin{subfigure}[b]{.32\linewidth}
				\centering
				\begin{tikzpicture}[scale=0.6,
	back/.style={dotted, thin},
	edge/.style={color=black},
	facet/.style={fill=Emerald,fill opacity=0.200000},
	vertex/.style={},
	back2/.style={dotted},
	edge2/.style={color=red,dashed},
	facet2/.style={fill=Emerald,fill opacity=0.20000},
	vertex2/.style={inner sep=1pt,circle,draw=darkgray!25!black,fill=darkgray!75!black,thick},
	back3/.style={dotted, thin},
	edge3/.style={color=black,very thick},
	facet3/.style={fill=Emerald,fill opacity=0.10000},
	vertex3/.style={inner sep=1pt,circle,draw=darkgray!25!black},
	facet4/.style={fill=blue,fill opacity=0.20000},]
	\def\x{2.5};
	\def\y{0};
	
	\draw[back] (-1, -\x-\y-1) to (1,-\x-\y-1) to (1,\x+\y+1) to (-1,\x+\y+1) to (-1, -\x-\y-1);
	
	\draw[back] (-\x-\y-1, 1) to (-\x-\y-1,-1) to (\x+\y+1,-1) to (\x+\y+1,1) to (-\x-\y-1, 1);

	\fill[facet] (-\x-1,-\y+1) to (-\x-1,-\y-1) to (-\x+1,-\y-1) to (\x+1,\y-1) to (\x+1,\y+1) to (\x-1, \y+1) to (-\x-1,-\y+1);
	\draw[edge] (-\x-1,-\y+1) to (-\x-1,-\y-1) to (-\x+1,-\y-1) to (\x+1,\y-1) to (\x+1,\y+1) to (\x-1, \y+1) to (-\x-1,-\y+1);
	
	\draw[edge2] (-1, -\x-\y-1) to (-\x-\y-1,-1);
	\draw[edge2] (1, -\x-\y-1) to (-\x-\y-1,+1);
	\draw[edge2] (1, \x+\y+1) to (\x+\y+1,1);
	\draw[edge2] (-1, \x+\y+1) to (\x+\y+1,-1);
	
	\node[vertex2] at (\x,\y) {};
	\node at (\x+0.4,\y+0.4) {$y$};
	\node[vertex2] at (-\x,-\y) {};
	\node at (-\x-0.4,-\y-0.4) {$x$};
\end{tikzpicture}
				
				\caption{end of movement $Z_{\max I}$}
			\end{subfigure}
			\caption{The SRS-movement $(Z_t)_{i \in I}$. The points $x$ and $y$ move along the red dashed line segments.}
		\end{figure}
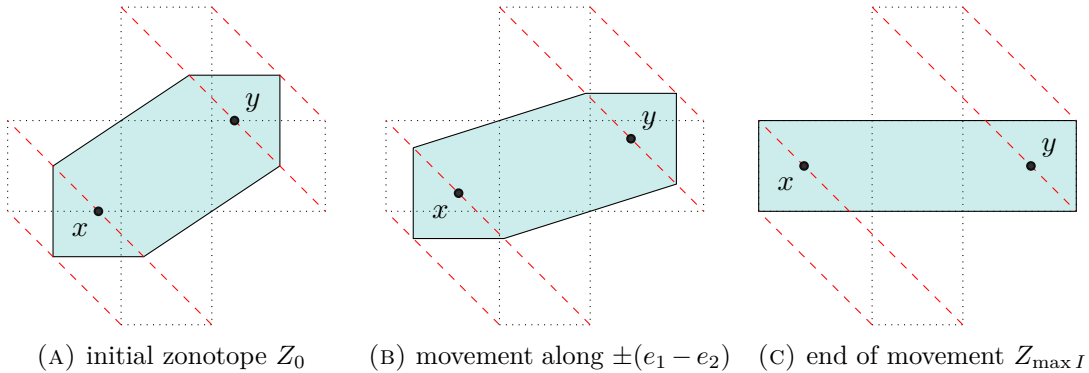
	
		The families $(Z_t \cap P_-)_t$ and $(Z_t \cap P_+)_t$ are RS-movements with common direction $[1, -1]^\transp$ and constant speed functions equal to $s_x$ and $s_y$, respectively. The remaining part $Z_t \cap P_0$ is a parallelogram with two opposing sides moving with speeds $s_x$ and $s_y$ along the lines $f^{-1}(\{x_1+x_2\})$ and $f^{-1}(\{y_1+y_2\})$, respectively. Therefore, if $x \neq y$ or $s_x = s_y$, the unique continuous speed function $\beta$ which is equal to $s_x$ on $P_-$, equal to $s_y$ on $P_+$ and affine on $P_0$ describes $Z_t$ via the formula \eqref{eq_def_rs_movement}. If $x = y$ and $s_x \neq s_y$, then $I = \{0\}$ and $(Z_t)_{t\in I}$ is an RS-movement for any speed function $\beta$.
	\end{proof}

	Building on Lemma \ref{lemma_2d}, we now address all dimensions $n \geq 2$.
	
	\begin{lemma} \label{lemma_higher_n}
		For $n \geq 2$, let $z \in \R^n$ with $z\geq 0$. We fix two indices $i < j \in [n]$ define $v_{ij}\in \R^n$ by $v_{ij} \coloneqq e_i-e_j$. For \begin{equation}
			t \in I \coloneqq \{r \in \R \mid z+rv_{ij} \geq 0\},
		\end{equation} 
		we define
		\begin{equation}
			Z_t \coloneqq C_n + [-(z+tv_{ij}),z+tv_{ij}].
		\end{equation}
		Then $(Z_t)_{t \in I}$ is an RS-movement. If $\max(z_i,z_j) >0$, then the corresponding speed function is not affine.
	\end{lemma}
	
	\begin{proof}
		By definition, $Z_t$ is a zonotope for all $z \in I$. As in the proof of Lemma \ref{lemma_2d}, we have to show that there is a speed function $\beta$ such that $Z_t$ can be written as in \eqref{eq_def_rs_movement}.	Without loss of generality, we assume that $i=1$ and $j=2$. %
		We fix $t\in I$ and $w \in \linh\{e_3,\dots,e_n\}$ and consider the polytope
		\begin{equation}
			P_t(w)\coloneqq Z_t \cap (w+\linh\{e_1,e_2\}),
		\end{equation}
		which is given by
		\begin{align}
			C_2& + \left[\left([-(z+tv_{12}),z+tv_{12}]+\sum_{i=3}^n[-e_i,e_i]\right) \cap (w+\linh\{e_1,e_2\}) \right] \\
			=C_2&+ \pi_w\left[[-(z+tv_{12}),z+tv_{12}] \cap \left(w+\linh\{e_1,e_2\}+\sum_{i=3}^n[-e_i,e_i]\right) \right],
		\end{align}
		where $C_2 = \sum_{i=1}^2[-e_i,e_i] \subset \R^n$ and $\pi_w$ denotes the orthogonal projection onto $w+\linh\{e_1,e_2\}$. If the second Minkowski summand in the right-hand side is empty, then $P_t(w)=\emptyset$. Otherwise, there exist $-1 \leq s_x^w \leq s_y^w \leq 1$ such that
		\begin{align}
			P_t(w)&=\sum_{i=1}^2[-e_i,e_i]+ \pi_w\left([s_x^w(z+tv_{12}),s_y^w(z+tv_{12})] \right)\\
			&=\sum_{i=1}^2[-e_i,e_i]+ \left[\pi_w(s_x^wz)+ts_x^wv_{12},\pi_w(s_y^wz)+ts_y^wv_{12} \right].
		\end{align} 
		Since $v_{12}\in \linh\{e_1,e_2\}$, the numbers $s_x^w$ and $s_y^w$ are independent of $t$. Identifying $w+\linh\{e_1,e_2\}$ with $\R^2$, we are precisely in the situation of Lemma \ref{lemma_2d} and obtain the existence of a speed function $\beta_w \colon w+\linh\{e_1,e_2\} \rightarrow \R$ such that $P_t(w)$ is an RS-movement of $P_0(w)$ on the interval
		\begin{equation}
			\tilde{I}\coloneqq\{t \in \R \mid \pi_w(s_x^wz)+ts_x^wv_{12} \leq \pi_w(s_y^wz)+ts_y^wv_{12}\}.
		\end{equation}
		If $s_x^w = s_y^w$, then clearly $\tilde{I}=\R$; otherwise a simple computation shows that $\tilde{I}=I$. Combining all speed functions $\beta_w$ for $w \in \linh\{e_3,\dots,e_n\}$ into a single speed function $\beta \colon \R^n \rightarrow \R$ via $\beta(w+x)\coloneqq\beta_w(x)$, $x \in \linh\{e_1,e_2\}$, we obtain that $(Z_t)_{t \in I}$ is an RS-movement. Finally, if $\max(z_1,z_2)>0$, then $I = [a,b] \neq \{0\}$. Since the zonotopes $Z_a$ and $Z_0$ are not combinatorially equivalent, they are in particular not affinely equivalent and $\beta$ cannot be affine.
	\end{proof}
	
	The proof of Proposition \ref{prop_n+1_zonotopes_upper_bound} is now merely a matter of combining the previous results.
	
	\begin{proof} [Proof of Proposition \ref{prop_n+1_zonotopes_upper_bound}]
		Since $Z$ is full-dimensional, it has $n$ linearly independent generators and we can assume without loss of generality that
		\begin{equation}
			Z = C_n + [-z,z]
		\end{equation}
		for some $z \geq 0$. We first consider the case $n=2$. If $Z$ is not affinely equivalent to $C_2$, then $z>0$. Setting $s_x \coloneqq 1$, $s_y \coloneqq -1$, $x \coloneqq -z$ and $y \coloneqq z$, Lemma \ref{lemma_2d} yields an RS-movement $(Z_t)_{t \in [a,b]}$ where both $Z_a$ and $Z_b$ are rectangles. Because the corresponding speed function is not affine, $t \mapsto L_{Z_t}^{4}$ is strictly convex, which implies $L_{C_2}=L_{Z_a}>L_{Z_0}=L_Z$.
		
		Now let $n \geq 3$. If $Z$ is not affinely equivalent to $C_n$, there exist two indices $i \neq j$ with $z_i>0$ and $z_j>0$. Then Lemma \ref{lemma_higher_n} yields the existence of an RS-movement $(Z_t)_{t \in [a,b]}$ with a non-affine speed function, where $Z_b$ is given by $C_n + [-\tilde{z},\tilde{z}]$ with $\tilde{z}_i=z_i+z_j$ and $\tilde{z}_j=0$. The other endpoint $Z_a$ of the RS-movement is affinely equivalent to $Z_b$ and is obtained by permuting the $i$-th and $j$-th coordinates of $Z_b$. 
		Because the corresponding speed function is not affine, $t \mapsto L_{Z_t}^{2n}$ is strictly convex, implying $L_{Z_a}=L_{Z_b}>L_{Z_0}=L_Z$. Iterating this procedure on possible further pairs of non-zero coordinates of $\tilde{z}$, we eventually arrive at a cuboid and obtain that  $L_{C_n}>L_Z$.
	\end{proof}

	In the proof of Proposition \ref{prop_n+1_zonotopes_upper_bound}, we iteratively transformed the zonotope $Z$ into a cuboid while strictly increasing the isotropic constant in each step. Clearly, already the first step suffices to show that the initial zonotope $Z$ cannot be a local maximizer. To deduce Proposition \ref{prop_n+1_zonotopes_upper_bound}, one could alternatively use the following lemma, which asserts that a global maximizer of $K \mapsto L_K$ in the class of zonotopes with at most $n+1$ generators indeed exists. Lemma \ref{lemma_higher_n} then shows that any zonotope with $n+1$ irredundant generators cannot be this maximizer, and Proposition \ref{prop_n+1_zonotopes_upper_bound} follows.
	
	\begin{lemma} \label{lemma_n+1_zonotopes_compactness}
		Let $Z \subset \R^n$ be a zonotope with at most $n+1$ generators. Then there exists a vector $y \in [0,1]^n$ such that
		\begin{equation}
			\tilde{Z} \coloneqq C_n +[-y,y]
		\end{equation}
		is affinely equivalent to $Z$. In particular, the function
		\begin{equation}
			\{Z \subset \R^n \mid Z \text{ is a zonotope with at most $n+1$ generators}\} \rightarrow \R, \quad Z \mapsto L_Z
		\end{equation}
		attains its minimum and its maximum.
	\end{lemma}
	
	\begin{proof}
		Let $Z$ be given by $Z = c_Z+\sum_{i=1}^{n+1} [-z_i,z_i]$ for some possibly redundant generators $z_1,\dots,z_n \in \R^n$ and a point $c_Z \in \R^n$. Since $\dim Z = n$, the matrix $M \coloneqq [z_1,\dots,z_{n+1}]$ has rank $n$ and there exists a linear dependence $0 = \sum_{i=1}^{n+1} \lambda_i z_i$ with at least one non-zero coefficient $\lambda_i$. Because we can permute the indices and replace $z_i$ by $-z_i$ if necessary, we can assume without loss of generality that $\lambda_{n+1} \geq \lambda_i \geq 0$ for all $i\in [n]$. Then we have $z_{n+1} = \sum_{i=1}^{n}\frac{\lambda_i}{\lambda_{n+1}} z_i$
		and the affine bijection that maps $0$ to $c_Z$ and $e_i$ to $c_Z+z_i$, $i \in [n]$, is an affine isomorphism from
		\begin{equation}
			\tilde{Z} =C_n +\left[-\sum_{i=1}^{n}\frac{\lambda_i}{\lambda_{n+1}}e_i,\sum_{i=1}^{n}\frac{\lambda_i}{\lambda_{n+1}}e_i\right]
		\end{equation}
		to $Z$.
		
		The second claim of the lemma follows from the compactness of $[0,1]^n$ and the continuity of $y \mapsto L_{C_n+[-y,y]}$ in combination with the fact that $L_Z$ is an affine invariant.
	\end{proof}
	
	Turning to the minimizers of $K \mapsto L_K$ in the class of zonotopes with at most $n+1$ generators, we consider the following infinite family of zonotopes:
	\begin{equation}
		Q_n \coloneqq C_n + [-(e_1+e_2+\dots+e_n),e_1+e_2+\dots+e_n]\subset \R^n \quad \text{for } n \geq 2.
	\end{equation}
	We note that $Q_2$ is affinely equivalent to a regular hexagon, whereas $Q_3$ is affinely equivalent to a rhombic dodecahedron. These highly symmetric polytopes and their $n$-dimensional relatives turn out to be the minimizers that we are looking for.
	\begin{proposition} \label{prop_n+1_zonotopes_lower_bound}
		Let $Z \subset \R^n$ be a zonotope with at most $n+1$ generators. Then
		\begin{equation}
			L_Z \geq L_{Q_n}
		\end{equation}
		with equality if and only if $Z$ is affinely equivalent to $Q_n$.
	\end{proposition}

	\begin{proof}
		Let $Z \subset \R^n$ be a global minimizer of $K \mapsto L_K$ in the class of zonotopes with at most $n+1$ generators, which exists by Lemma \ref{lemma_n+1_zonotopes_compactness}. Without loss of generality, let $Z$ be given by $Z=C_n +[-y,y]$ for a vector $y \geq 0$. We assume towards a contradiction that here exist indices $i<j$ with $y_i \neq y_j$. In particular, we have $\max(y_i,y_j)>0$ and Lemma \ref{lemma_higher_n} yields an RS-movement $(Z_t)_{t \in [a,b]}$ of $Z\coloneqq C_n+[-y,y]$ with non-affine speed function. Therefore, by Theorem \ref{thm_ccg}, the function $t \mapsto L_{Z_t}^{2n}$ is strictly convex. Moreover, it is easy to see that the RS-movement $(Z_t)_{t \in [a,b]}$ is symmetric in the sense that $Z_{a+r}$ is affinely equivalent to $Z_{b-r}$ for all $r \in [0, b-a]$. Therefore, the function $t \mapsto L_{Z_t}$ attains its minimum precisely at the point $\frac{a+b}{2}$. Because $y_i \neq y_j$, we have $\frac{a+b}{2} \neq 0$ and it follows that $Z$ is not a minimizer of $K \mapsto L_K$, contradicting our assumption. It follows that $y_1=y_2=\dots=y_n=\lambda$ or, in other words,
		\begin{equation}
			y = \lambda(e_1+\dots+e_n) 
		\end{equation}
		for some $\lambda >0$.
		
		It remains to show that $\lambda=1$. From the proof of Proposition \ref{prop_n+1_zonotopes_upper_bound}, we know that $y \neq 0$. Therefore, $y,e_2,\dots,e_n$ forms a basis of $\R^n$ and we can repeat our construction of the RS-movement $(Z_t)_{t \in [a,b]}$ above, with $-e_1$ assuming the role of $y$ and $-y$ assuming the role of $e_1$. We obtain
		\begin{equation}
			-e_1 = \mu(- y +e_2 +\cdots + e_n).
		\end{equation}
		for some $\mu >0$, which implies $\lambda=y_n=1$.
	\end{proof}
	
	\subsection*{Acknowledgments}
	Funded by the Deutsche Forschungsgemeinschaft (DFG, German Research Foundation) under Germany's Excellence Strategy – The Berlin Mathematics Research Center MATH+ (EXC-2046/1, project ID: 390685689). I thank the anonymous reviewers, Tom Baumbach, Katharina Eller and Emanuel Milman for helpful comments and suggestions.

	\printbibliography
\end{document}